\documentclass[12pt]{article}
\usepackage{graphicx} 
\usepackage{fullpage}
\usepackage{tikz}
\usepackage{xurl}
\title{Linear Systems and Eigenvalue Problems: \\Open Questions from a Simons Workshop}
\author{
Noah Amsel${}^{\dagger}$ \and Yves Baumann${}^{\dagger}$ \and Paul Beckman${}^{\dagger}$ \and Peter B{\"u}rgisser${}^{\dagger}$ \and Chris Camaño${}^{\dagger}$ \and Tyler Chen${}^{\dagger}$ \and Edmond Chow${}^{\dagger}$ \and Anil Damle${}^{\dagger}$ \and Michal Derezinski${}^{\dagger}$ \and Mark Embree${}^{\dagger,\ast}$ \and Ethan N.\ Epperly${}^{\dagger}$ \and Robert Falgout${}^{\dagger}$ \and Mark Fornace${}^{\dagger}$ \and Anne Greenbaum${}^{\dagger}$ \and Chen Greif${}^{\dagger}$ \and Diana Halikias${}^{\dagger}$ \and Zhen Huang${}^{\dagger}$ \and Elias Jarlebring${}^{\dagger}$ \and Yiannis Koutis${}^{\dagger}$ \and Daniel Kressner${}^{\circ}$ \and Rasmus Kyng${}^{\dagger,\diamond}$ \and J\"org Liesen${}^{\dagger,\ast,\diamond}$ \and Jackie Lok${}^{\dagger}$ \and Raphael A.\ Meyer${}^{\dagger,\ast}$ \and Yuji Nakatsukasa$^{\ast}$ \and Kate Pearce${}^{\dagger}$ \and Richard Peng${}^{\dagger}$ \and David Persson${}^{\dagger}$ \and Eliza Rebrova${}^{\dagger}$ \and Ryan Schneider${}^{\dagger}$ \and Rikhav Shah${}^{\dagger}$ \and Edgar Solomonik${}^{\dagger,\ast,\diamond}$ \and Nikhil Srivastava${}^{\dagger,\diamond}$ \and Alex Townsend${}^{\dagger}$\and Robert J. Webber${}^{\dagger}$ \and Jess Williams${}^{\dagger}$ \footnote{In the list of authors, $\dagger$ denotes scribe, $\ast$ denotes working group lead, $\diamond$ denotes editor, and $\circ$ denotes lead organizer.}
}

\newcommand{\scribe}[1]{\noindent \emph{Scribe: #1}\\}
\newcommand{\source}[1]{\noindent \emph{Source: #1}\\}
\newcommand{\motivation}{\emph{Motivation.} }
\newcommand{\remarks}{\paragraph{Remarks.}}
\newcommand{\poly}{\mathrm{poly}}

\newenvironment{update}
  {\par\noindent\color{red}Update on August 20, 2026: \ignorespaces}
  {\par}

\usepackage{hyperref}
\usepackage{amsmath,amsthm,amssymb,mathtools}
\usepackage{mleftright}
\newtheorem{problem}{Problem}[section]
\newtheorem{theorem}{Theorem}[section]

\newtheorem*{theorem*}{Theorem}

\newtheorem*{lemma*}{Lemma}

\newtheorem*{prop*}{Proposition}

\newcommand{\C}{{\mathbb{C}}}
\newcommand{\R}{{\mathbb{R}}}

\newcommand{\E}{{\mathbb{E}}}
\newcommand{\K}{{\mathcal{K}}}

\newcommand{\eps}{\varepsilon}
\newcommand{\abs}[1]{\mleft|#1\mright|}

\DeclareMathOperator*{\area}{area}

\makeatletter
\renewcommand{\@fnsymbol}[1]{\@arabic{#1}}
\makeatother

\theoremstyle{definition}
\newtheorem{definition}[theorem]{Definition}
\usepackage{soul,xcolor}

\begin{document}
\maketitle
\setcounter{footnote}{1}

\begin{abstract}
This document presents a series of open questions arising in matrix computations, i.e., the numerical solution of linear algebra problems. It is a result of working groups at the workshop \emph{Linear Systems and Eigenvalue Problems}, which was organized at the Simons Institute for the Theory of Computing program on \emph{Complexity and Linear Algebra} in Fall 2025. The complexity and numerical solution of linear algebra problems is a crosscutting area between theoretical computer science and numerical analysis. The value of the particular problem formulations here is that they were produced via discussions between researchers from both groups.

The open questions are organized in five categories: iterative solvers for linear systems, eigenvalue computation, low-rank approximation, randomized sketching, and other areas including tensors, quantum systems, and matrix functions.

\textcolor{red}{Updated to reflect the status of the open problems as of August 20, 2026.}
   
\end{abstract}

\maketitle

\newpage

\tableofcontents

\newpage

\section{Introduction}

Numerical algorithms are a central tool for modeling physical phenomena as well as in machine learning (data science).
Development of the theory and practice for numerical problems is a challenge that involves design of exact or approximate algorithms with a bounded number of steps, as well as considerations of the effects of finite precision and bit complexity.
The algorithmic landscape for basic numerical linear algebra problems such as solving linear systems and eigenvalue problems has evolved significantly in the last decade largely due to increasing prevalence of approximation by randomized mechanisms such as matrix sketching.

We organize proposed open questions (and associated discussion of context and general research directions) in sections that correspond to independent working groups 
at the workshop \emph{Linear Systems and Eigenvalue Problems}, held October 6-10, 2025 at the Simons Institute for the Theory of Computing \footnote{\url{https://simons.berkeley.edu/workshops/linear-systems-eigenvalue-problems}}.
These discussions built upon a set of $55$ open questions proposed by participants ahead of and at the start of the workshop, and had the general goal of assembling a set of well-defined open questions with agreed-upon formulations. Each working group aimed to form discussions among sets of 5-10 people from a diverse set of research areas.

The first section of open problems (Section~\ref{sec:linsys}) concerns the design and analysis of iterative solvers for linear systems.
Iterative linear solvers are the most widely used numerical tools in solution of large-scale simulation problems such as those arising from discretization of partial differential equations describing fluid flow.
While extensively studied in numerical analysis, there remain many gaps in our understanding of their properties. Section~\ref{sec:linsys} puts forward research directions regarding the development of more advanced representative benchmark problems (in terms of discretizations and matrix properties) for this purpose.
A number of problems concern the numerical behavior and convergence properties of the multigrid method, conjugate gradient, and the generalized minimum residual method (GMRES). 

Section~\ref{sec:eig} considers eigenvalue problems. In this area, randomized algorithms have recently shown to be effective in providing a provably convergent variant of QR iteration for general matrices, a long standing open problem in numerical analysis. We consider the problem of derandomizing this approach (and more generally designing deterministic matrix perturbations that achieve similar spread of the eigenvalue spectrum as the randomized case). Additionally, we survey open challenges in theoretical analysis related to stable and efficient computations of the eigenvalues of tridiagonal matrices and singular values of bidiagonal matrices. In the context of iterative solvers, we also pose questions related to understanding of the convergence of Krylov methods (distribution of Ritz values).

Low-rank approximation is another pillar of numerical linear algebra, and questions on this topic are covered in Section~\ref{sec:lowrank}.
We consider open problems regarding understanding of the effectiveness of standard (greedy) methods, such as QR with column pivoting or LU with complete pivoting, as well as randomized schemes.
Again, we also give focus to classes of matrices that arise from standard continuous solvers or satisfy common structural properties (including Laplacians and diagonally dominant matrices). A number of questions specifically concern low-rank approximation by column subset selection.
Accuracy properties of algorithms for hierarchical and block matrix approximations are considered here, and also later revisited for the tensor case in Section~\ref{sec:other}.

Section~\ref{sec:sketch} studies the properties of randomized embeddings and their application in numerical linear algebra (including those of focus in prior sections).
A recent development in the theory of subspace embeddings is the study of injections, which satisfy looser properties than the standard notion of randomized subspace embeddings.
Injections hence allow for reduced sample/sketch size, and yet are sufficient for effective approximations for important problems such as linear least squares.
We pose problems regarding their effectiveness in a number of other problems.
Section~\ref{sec:sketch} also considers the problem of optimality of known sparse embeddings, and a more complete characterization of the properties of structured embeddings (in particular subsampled randomized Hadamard transforms).

Section~\ref{sec:other} considers a few problems beyond core topics of classical and randomized numerical linear algebra.
First, we consider the problem of devising general tensor decomposition algorithms with provably bounded approximation error.
Additionally, we consider a model problem for eigenvector-dependent nonlinear eigenvalue problems (a class of optimization problems that forms the basis of widely used methods in quantum chemistry) based on the 2-local Hamiltonian eigenvalue problem.
Finally, we consider the characterization of a class of matrix functions that correspond to maps that can be computed by a fixed number of matrix-matrix products (hence of practical importance), and the approximation properties of this class of functions.

\section{Linear algebraic systems}
\label{sec:linsys}

\subsection{Scalable iterative solvers for a set of parameterized benchmark problems}\label{sec:scalable}

\scribe{Yves Baumann, Edmond Chow, Robert Falgout, Yiannis Koutis, and Jörg Liesen}

\emph{Motivation.}
The theoretical computer science (TCS) community and the numerical
linear algebra (NLA) community both work on algorithms for solving
linear systems. A closer collaboration between these communities
could be fueled by defining a set of benchmark linear system problems
that are of interest to the NLA community and that will expand the
range of problems for the TCS community to address.

Theoretical computer science has developed algorithms for solving
certain classes of sparse linear systems of equations that build
on graph Laplacian solvers.  These classes include symmetric
diagonally dominant M-matrices (SDDM matrices) and more general
symmetric diagonally dominant (SDD) matrices. In the latter, some
off-diagonal entries can be positive. An important next step in
this development is to extend the classes of sparse linear systems
that are scalably solvable (near linear convergence rates independent
of problem size) with the techniques of theoretical computer science.

We propose focusing the benchmark problems on linear systems that
are discretizations of partial differential
equations (PDEs) due to both their wide application in science and
engineering and to their connection to graph Laplacians that have
been a focus of TCS research in linear system solvers.
Furthermore, matrix collections such as the SuiteSparse collection
\cite{DavisHu2011} provide a set of matrices, but they do not
provide a pipeline that defines parametrized discretizations of
model problems from which a large set of such sample matrices could
be generated. Parameterized classes of problems are important to
understand the range of characteristics and behaviors that a solver
might encounter for a class of problem.

We therefore propose the following two problems.

\begin{problem}
Develop software to generate benchmark parameterized linear systems
that are generally accepted as relevant in practice and suitable
for theoretical analysis.
\end{problem}

\begin{problem}
For each class of parameterized benchmark linear systems, develop
theoretically and practically fast algorithms.
\end{problem}

Below, we define some parameterized PDEs that can lead to parameterized
classes of matrix discretizations to serve our purpose.  These
problems were chosen because we believe they are within reach for
algorithms with linear or nearly linear convergence
guarantees, possibly by extending techniques already developed in
TCS, and perhaps incorporating ideas developed in NLA.  Some of
these problems could also be challenging to solve by NLA techniques
if an entire matrix class is to be solved {\em robustly} and/or in
{\em linear time} for a very wide range of problem parameter values.
However, these problems are not intended to represent the most
challenging problems for PDE solvers, which might involve large
systems of PDEs modeling different physics on different scales.

\medskip
\emph{Diffusion problems.}
The diffusion equation is given by
\begin{equation}
    -\nabla \cdot (a(x,y) \nabla u) = f(x,y)  \text{ in } \Omega,
\end{equation}
where $u$ is an unknown scalar function, $a$ is tensor function of
diffusion coefficients, $f$ is a scalar function of sources and sinks, and
$\Omega$ defines a domain.
Various boundary conditions are also needed to specify the problem,
but zero Dirichlet conditions could be considered for simplicity.
Finite difference and finite element discretizations of this
PDE lead to linear systems of equations.

The problem is parameterized by the diffusion tensor $a$ which can
make the problem more or less challenging to solve.  One challenging
family of diffusion problems arises when $a$ has large jumps.
Specific example problems from NLA could be found in
\cite{Grigori2014,Grigori2021}, among others.  Furthermore, higher order
discretizations (beyond second order) of the diffusion PDE lead to
discrete matrices with {\em positive} off-diagonal elements and
other properties that could make them challenging for existing
graph Laplacian solvers.

Reference solvers from the NLA community for diffusion problems
include the BoomerAMG multigrid method \cite{Falgout2002} and
the GENEO preconditioner with its implementation in the HPDDM library
\cite{Jolivet2013}.

\medskip
\emph{Convection-diffusion problems.}
Convection-diffusion PDEs give rise to nonsymmetric systems of equations.  Discretized convection-diffusion equations are also often subproblems for more complex PDEs, such as the Navier-Stokes equations.
The convection-diffusion PDE is
\[
-\nabla \cdot (a(x,y) \nabla u) + b(x,y) \cdot \nabla u = f(x,y),
\]
where $u$ is an unknown scalar function, $a$ is tensor function of diffusion coefficients, $b$ is a vector function of convection coefficients, and $f$ is a scalar function of sources and sinks. A domain and various boundary conditions are also needed to specify the problem, but a square domain and zero Dirichlet conditions could be considered for simplicity.

The PDE can be discretized on a regular mesh using finite differences.  For irregular geometries and irregular meshes, more sophisticated discretizations such as finite elements must be used.

For simplicity, consider a finite difference discretization of a constant coefficient convection-diffusion equation on a 2-D regular mesh with $a= \left[ \begin{smallmatrix} a_x & 0 \\ 0 & a_y \end{smallmatrix} \right]$ (which must be positive semidefinite) and $b=[b_x,b_y]^T$. Using a second order approximation of the second derivatives and first order centered approximation for the first derivatives, we obtain the following 5-point stencil operator at each grid point, assuming that the mesh spacing $h=1$.
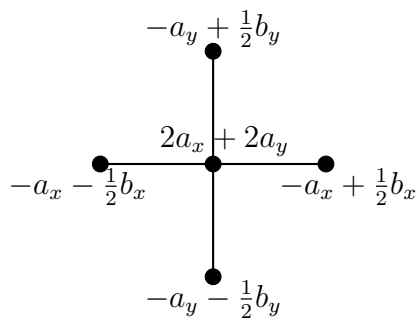
\begin{center}
\begin{tikzpicture}[scale=1.5]
  \filldraw[black] (0,0) circle (0.07); 
  \filldraw[black] (-1,0) circle (0.07); 
  \filldraw[black] (1,0) circle (0.07); 
  \filldraw[black] (0,-1) circle (0.07); 
  \filldraw[black] (0,1) circle (0.07); 

  \draw[thick] (0,0) -- (-1,0); 
  \draw[thick] (0,0) -- (1,0);  
  \draw[thick] (0,0) -- (0,-1); 
  \draw[thick] (0,0) -- (0,1);  

  \node at (0.1,0.2)  {$2 a_x + 2 a_y$};
  \node at (-1.2,-0.2) {$-a_x -\tfrac{1}{2}b_x$};
  \node at ( 1.2,-0.2) {$-a_x +\tfrac{1}{2}b_x$};
  \node at (0,-1.2) {$-a_y -\tfrac{1}{2}b_y$};
  \node at (0,1.2)  {$-a_y +\tfrac{1}{2}b_y$};
\end{tikzpicture}
\end{center}
Notice that the convection terms give rise to nonsymmetry in the stencil.  Each interior mesh point corresponds to a matrix row that maintains weak diagonal dominance when the convection coefficients are small (and all off-diagonal entries remain nonpositive).

When the convection coefficients are large in magnitude, the off-diagonal entries can become positive and diagonal dominance is lost.  In this case, the discretization is unstable.  This can be fixed by using one-sided upwind differences to approximate the first derivatives.  The choice of which side of the mesh point to take the one-sided difference depends on the sign of $b_x$ or $b_y$.  In short, the difference is chosen so that a positive contribution is made to the diagonal of the matrix.  Assuming $b_x$ and $b_y$ are positive, the stencil when using upwind differencing is the following.
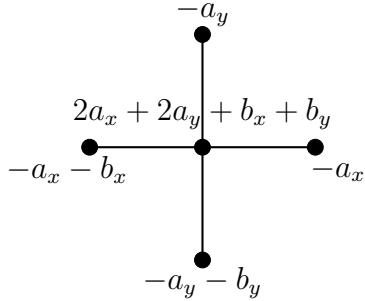
\begin{center}
\begin{tikzpicture}[scale=1.5]
  \filldraw[black] (0,0) circle (0.07); 
  \filldraw[black] (-1,0) circle (0.07); 
  \filldraw[black] (1,0) circle (0.07); 
  \filldraw[black] (0,-1) circle (0.07); 
  \filldraw[black] (0,1) circle (0.07); 

  \draw[thick] (0,0) -- (-1,0); 
  \draw[thick] (0,0) -- (1,0);  
  \draw[thick] (0,0) -- (0,-1); 
  \draw[thick] (0,0) -- (0,1);  

  \node at (0.0,0.3)  {$2 a_x + 2 a_y + b_x + b_y$};
  \node at (-1.2,-0.2) {$-a_x -b_x$};
  \node at ( 1.2,-0.2) {$-a_x$};
  \node at (0,-1.2) {$-a_y -b_y$};
  \node at (0,1.2)  {$-a_y$};
\end{tikzpicture}
\end{center}

Notice that the matrix remains diagonally dominant.

More challenging examples of convection-diffusion problems do not necessarily fix $b$ and may also use unstructured meshes for their discretization.  Example problems are available from \cite[Chapter~6]{ElmanHowardSilvester2014}.  The example problems here set $a$ to a constant scalar $\epsilon$ and specify the advection field $b$.  Possible examples include $b=(0,1)^T$ in 2-dimensions, which has an analytical solution under appropriate boundary conditions; a 2-D field that varies quadratically in the vertical dimension; a constant advection field that is not aligned with the grid; and a recirculating advection field.  Reference solvers are available in the associated library IFISS \cite{ifiss2016}.

For nonsymmetric linear equations such as those arising from convection-diffusion equations, efficient short-recurrence Krylov subspace methods such as the conjugate gradient method are not available.  The convergence of methods such as GMRES for nonsymmetric equations depend not only on the eigenvalues of the (preconditioned) matrix, but also on the conditioning of the eigenvectors. Thus TCS techniques that are based on Krylov subspace methods may not be amenable for nonsymmetric equations. Instead, approaches using stationary iterations may be more susceptable to proof. In all cases, the desired accuracy of the approximate solution should be comparable to that of the discretization error, $O(h^2)$, measured in the same norm.

\medskip
\emph{Helmholtz equation.}
The continuous Helmholtz equation is
\[
-\nabla^2 u - k^2 u = f,
\]
where $u$ is an unknown scalar function, and $k$ is called the wave number.
For simplicity for the purposes of this document, we may consider Dirichlet
boundary conditions on the domain.  Discretization of this PDE, for example,
by finite-differences, leads to a matrix equation of the form
\[
(A - k^2 I) x = b,
\]
where $A$ is a symmetric positive definite matrix representing the discretization
of $-\nabla^2 u$.  One may need to solve this equation for several values of $k$ such that the matrix $(A-k^2 I)$ has several negative eigenvalues.  Thus the discrete problem is indefinite, which is considered to be challenging for iterative solution methods \cite{ErnstGander2011}.

A well-developed approach for solving the discrete equations efficiently is to use a {\em shifted Laplacian preconditioner} combined with an iterative method.  The shifted Laplacian is
\[
A - (\alpha + i \beta) k^2 I,
\]
where $\alpha$ and $\beta$ are real parameters that help shift the spectrum of the preconditioned matrix away from zero.  Note that the shift is complex. Solving with the shifted Laplacian can be done by solving with its incomplete factors \cite{Erlangga2004}.

\subsection{Correctness of multigrid methods beyond the standard model problems}

\scribe{Edmond Chow, Robert Falgout, and Rasmus Kyng} 

\emph{Motivation.}
In NLA, multigrid methods are optimal $O(n)$ solvers for $n$ degrees of freedom (see, e.g., \cite{Briggs2000,Hac85,Trottenberg2001} for general introductions to multigrid methods).
Optimal $O(n)$ complexity is also the goal for solvers in TCS.
Convergence proofs from NLA are available for discretizations of various boundary value
problems using multigrid methods designed for these problems.  The conjecture we propose here is to prove the $O(n)$ complexity of a multigrid method for a general class of linear systems with polynomially bounded SWCDDM matrices, described algebraically below.  Due to their similarity to discretizations of elliptic PDEs, proving multigrid convergence for these problems may be possible.  For example, an approach may be to prove using TCS techniques that the smoothing and approximation properties (see below) are satisfied with properly chosen smoothers and prolongators chosen for a polynomially bounded SWCDDM matrix.

\begin{definition}
    \noindent
Consider a $\mathbf{M} \in R^{n \times n}$, and let $\mathbf{M} = \mathbf{D} - \mathbf{A}$
where $\mathbf{D}$ is diagonal and $\mathbf{A}$ is has zero diagonal.
    \begin{enumerate}
        \item $\mathbf{M}$ is symmetric diagonally dominant (SDD) if
         $\mathbf{M}$ is symmetric and for all $i$,
    $\mathbf{M}(i,i) \geq \sum_{j \neq i } |\mathbf{M}(i,j)|$.
        \item $\mathbf{M}$ is SDDM if it is SDD and has non-positive off-diagonals.
        \item A row $i$ of $\mathbf{M}$ is \emph{strictly} diagonally dominant if 
        $\mathbf{M}(i,i) > \sum_{j \neq i } |\mathbf{M}(i,j)|$ and $\alpha$-diagonally dominant if 
        $\mathbf{M}(i,i) \geq \alpha + \sum_{j \neq i } |\mathbf{M}(i,j)|$. 
        \item A row $i$ of $\mathbf{M}$ is \emph{weakly} diagonally dominant if 
        $\mathbf{M}(i,i) = \sum_{j \neq i } |\mathbf{M}(i,j)|$.
        \item A matrix $\mathbf{M} \in R^{n \times n}$ is symmmetric weakly-chained diagonally dominant M (SWCDDM)
        if (1) it is SDDM and (2) for every weakly-diagonally dominant row $i$,
        there exists a sequence of indices $i = i_0, i_1, \ldots, i_k$ s.t. for all $j \in \{ 1, \ldots, k\}$
        we have $\mathbf{A}(i_{j-1},i_{j}) \neq 0$ and row $i_k$ is strictly diagonally dominant.
    \end{enumerate}
\end{definition}

\begin{definition}
    \noindent
We say $\mathbf{M}$ is \emph{$c$-polynomially bounded} if $\max_{i,j} |M(i,j)| \leq n^c$
and $|M(i,j)| n^c$ is integral.
A collection of matrices is \emph{polynomially bounded} if there exists a constant $c$, s.t. every matrix in the collection \emph{$c$-polynomially bounded}.
\end{definition}

\begin{problem}[Algebraic Multigrid Analysis for Theoretical Computer Scientists]
    \label{prob:TCS_AMG}
Show that for polynomially bounded SWCDDM matrices $\mathbf{M}$,
Algebraic Multigrid runs in time $\tilde{O}(\operatorname{nnz}(\mathbf{M}))$ to compute an (implicitly represented) symmetric matrix $\mathbf{Z}$,
which can be applied in $O(\operatorname{nnz}(\mathbf{M}))$ time, and satisfies the condition that
$ \Omega(1) \mathbf{M}^{-1} \preceq \mathbf{Z} \preceq O(1) \mathbf{M}^{-1}$.  In particular, the constants are independent of the problem size $n$.
\end{problem}
We also want to highlight some variants of this problem:
\begin{enumerate}
    \item Solve Problem \ref{prob:TCS_AMG}, but only for \emph{planar} matrices.
    \item Solve Problem \ref{prob:TCS_AMG}, but spend even less time computing the (implicit) operator representation, e.g.
$O(\operatorname{nnz}(\mathbf{M}) \poly(\log\log(\operatorname{nnz}(\mathbf{M}))))$ or even $O(\operatorname{nnz}(\mathbf{M}))$ time.
    \item Solve Problem \ref{prob:TCS_AMG}, but spend slightly more time applying the operator $\mathbf{Z}$, say,
    \\ $O(\operatorname{nnz}(\mathbf{M}) \poly(\log\log(\operatorname{nnz}(\mathbf{M}))))$ time.
\end{enumerate}

A multigrid method consists of a smoother, prolongation (interpolation) and restriction operators, and a coarse grid correction operator. Each of these components is defined at each grid level of the multigrid hierarchy.  In NLA, proofs of multigrid convergence depend on conditions on the smoother and prolongation operator, called the {\em smoothing property} and the {\em approximation property}.  The smoother and prolongation operators are designed to satisfy these conditions depending on the partial differential equation (PDE) boundary value problem whose discrete equations are to be solved \cite{Yserentant1993,Xu2017}. As an example, the classical convergence theorem for multigrid has the following form \cite{Hac85, Wesseling1992}.

\emph{Two-grid convergence theorem for an elliptic model problem.}
Consider the model problem
\begin{equation*}
-\Delta u = f \quad \text{in } \Omega, \quad u|_{\partial \Omega} = 0,
\end{equation*}
discretized with finite-differences on a grid with mesh size $h$. Let $A$ denote the SPD stiffness matrix. Define a \textit{two-grid method} that satisfies:

\begin{enumerate}
  \item \textit{Smoothing Property}: The chosen smoother (e.g., Gauss--Seidel or Jacobi) has a smoothing property if there exist a constant $C_S$ and a function $\eta(\nu)$ independent of $h$ such that:
  \begin{equation*}
  \| A S^{\nu} \|_2 \leq C_S h^{-2} \eta(\nu), \quad \eta(\nu) \rightarrow 0 \ \mathrm{for}\ \nu \rightarrow \infty,
  \end{equation*}
  where $S$ is the smoother iteration matrix.

  \item \textit{Approximation Property}: The interpolation operator $P$ and the restriction operator $R$ satisfy the approximation property if there exists a constant $C_A$ independent of $h$ such that:
  \begin{equation*}
  \| A^{-1} - P (A^H)^{-1} R \|_2 \leq C_A h^{2},
  \end{equation*}
  where $A^H$ is the coarse grid matrix.
\end{enumerate}

Then the \textit{two-grid method} with $\nu$ pre-smoothing steps, defined by the iteration matrix
\begin{equation*}
Q = (A^{-1} - P (A^H)^{-1} R)(A S^{\nu})
\end{equation*}
satisfies:
\begin{equation*}
\| Q \|_2 \leq C_S C_A \eta(\nu) < 1, \quad \forall \nu \ge \tilde{\nu},
\end{equation*}
for $\tilde{\nu}$ independent of $h$. Thus, the convergence rate is uniform with respect to $h$.  This result can be extended to the multigrid V-cycle for multiple levels to show that the multigrid method for the model problem achieves optimal complexity $O(n)$ for $n$ unknowns.

\subsection{Conjugate gradient versus sketch-and-project} 
\scribe{Tyler Chen, Jackie Lok, and Robert J. Webber}

\emph{Motivation.}
Recently, researchers in randomized numerical linear algebra have developed many new methods for solving linear systems
\begin{equation*}
    \mathbf{A} \mathbf{x} = \mathbf{b},
    \qquad \text{where } \mathbf{A} \in \mathbb{C}^{n \times n} \text{ is strictly positive definite},
\end{equation*}
based on the sketch-and-project framework \cite{GowerRichtarik2015}.
This work raises the following question:
\begin{quote}
    What conditions on the eigenvalues of $\mathbf{A}$ ensure that a sketch-and-project method converges to $\varepsilon$ accuracy after a smaller number of epochs than a traditional Krylov solver, such as conjugate gradient \cite{HestenesStiefel1952}?
\end{quote}
One epoch of a sketch-and-project method is the time to read and process $n$ randomly selected rows.
Similarly, one epoch of a Krylov method is the time to perform a single matrix--vector multiplication, which reads and processes all $n$ rows of $\mathbf{A}$.

\emph{Scope.} 
Preconditioners can be used to speed up Krylov methods \cite{Saad2003}, just as they can be used to speed up sketch-and-project methods \cite{LokRebrova2025}.
Furthermore, sketch-and-project methods can be accelerated by processing many rows of $\mathbf{A}$ at a time \cite{derezinski2024sharp} and/or incorporating momentum \cite{GowerHaRiSt2018}.
The numerous bells and whistles raise interesting questions for analysis.
Yet we still need to understand the simplest possible comparison between the un-preconditioned conjugate gradient (CG) method and the randomized coordinate descent (RCD) method \cite{LeventhalLewis2010, Nesterov2012} defined as follows.

\begin{definition}[Randomized coordinate descent]
    Fix a strictly positive definite $n \times n$ linear system $\mathbf{A} \mathbf{x} = \mathbf{b}$.
    Starting with an initial iterate $\mathbf{x}_0$ and residual vector $\mathbf{r}_0 = \mathbf{b} - \mathbf{A} \mathbf{x}_0$, iterate over the following steps for $t = 1, 2, \ldots$:
    \begin{itemize}
        \item[(a)] Choose a random index $i_t \in \{1, \ldots, n\}$ with $\operatorname{prob}(i) = \mathbf{A}(i,i) / \operatorname{tr}(\mathbf{A})$.
        \item[(b)] Update the $i_t$\textsuperscript{th} entry of the current iterate so that $\mathbf{A}(i_t, \cdot) \mathbf{x}_t = \mathbf{b}(i_t)$ is satisfied exactly:
        \begin{equation*}
        \mathbf{x}_t = \mathbf{x}_{t-1} + \frac{\mathbf{r}_{t-1}(i_t)}{\mathbf{A}(i_t, i_t)} \mathbf{e}_{i_t}
        \quad \text{and} \quad
        \mathbf{r}_t = \mathbf{r}_{t-1} - \frac{\mathbf{r}_{t-1}(i_t)}{\mathbf{A}(i_t, i_t)} \mathbf{A}(\cdot, i_t).
        \end{equation*}
    \end{itemize}
\end{definition}

\emph{Problem statement.} 
Let $\mathbf{x}_{\star} = \mathbf{A}^{-1} \mathbf{b}$ be the solution to the linear system, and let $\lVert \cdot \rVert_{\mathbf{A}} = \langle \mathbf{A} \cdot, \cdot\rangle^{1/2}$ be the $\mathbf{A}$-weighted norm. 
Define the stopping times
\begin{align*}
    & T_{\rm CG}(\mathbf{A},\mathbf{b},\epsilon) = \min\{t \geq 0 \,{:}\, \text{CG achieves $\lVert \mathbf{x}_t - \mathbf{x}_{\star} \rVert_{\mathbf{A}}^2 \leq \varepsilon \lVert \mathbf{x}_0 - \mathbf{x}_{\star} \rVert_{\mathbf{A}}^2$ after $t$ steps} \}, \\
    & T_{\rm RCD}(\mathbf{A},\mathbf{b},\epsilon) = \min\{t \geq 0 \,{:}\, \text{RCD achieves $\lVert \mathbf{x}_t - \mathbf{x}_{\star} \rVert_{\mathbf{A}}^2 \leq \varepsilon \lVert \mathbf{x}_0 - \mathbf{x}_{\star} \rVert_{\mathbf{A}}^2$ after $t$ steps} \}.
\end{align*}
The following stylized problem is simple to state and it probes the difference between CG and RCD, but we have not been able to solve it:
\begin{problem}[Polynomially decaying eigenvalues]
\label{problem24}
    Fix constants $\varepsilon > 0$ and $p > 0$.
    Let $\mathbf{A}_n \in \mathbb{C}^{n \times n}$, $n\geq 1$, be a sequence of random positive definite matrices constructed via
    \begin{equation*}
        \mathbf{A}_n = \mathbf{U} \mathbf{\Lambda}_n \mathbf{U}^*,
        \quad \text{where } \begin{cases} \mathbf{U} \in \mathbb{C}^{n \times n} \text{ is Haar unitary}, \\
        \mathbf{\Lambda}_n = \operatorname{diag}\{1^{-p}, 2^{-p}, \ldots, n^{-p}\}.
        \end{cases}
    \end{equation*}
    Let $\mathbf{b}_n = \mathbf{A}_n \mathbf{z}_n$ be a sequence of random vectors where $\mathbf{z}_n \in \mathbb{C}^n$ is distributed uniformly on the sphere.
    What is the asymptotic behavior of $T_{\rm CG}(\mathbf{A}_n,\mathbf{b}_n,\epsilon)$ compared with $T_{\rm RCD}(\mathbf{A}_n,\mathbf{b}_n,\epsilon)$ as $n \rightarrow \infty$?
\end{problem}

\begin{update}
    The reference \begin{quote}
        Leheng Chen, Zihao Liu, Wanyi He and Bin Dong, Iteris: Agentic Research Loops for Computational Mathematics, arXiv:2606.02484, 2026
    \end{quote}
    claims to solve Problem~\ref{problem24}.

    As mentioned above, Problem~\ref{problem24} is a stylized question.  The much broader motivating question, characterizing when sketch-and-project beats CG for general spectra, right-hand sides, preconditioning and finite precision, remains open.
\end{update}

\subsection{The effect of outlying singular values (or eigenvalues)}

\scribe{Michal Derezinski} 

The performance of iterative methods  for solving a linear system $Ax=b$ (such as Krylov methods, or randomized sketching methods) often depends on the distribution of the eigenvalues of $A$, or its singular values if solved via the normal equations, beyond just the condition number of the matrix. A broad goal of this direction is:
\begin{quote} 
Find natural characterizations of the complexity of solving linear systems in terms of the properties of the eigenvalue or singular value distribution of the input.
\end{quote}

\emph{Outlying singular values.} A commonly occurring scenario is that of matrices with a spectrum that consists of a large cluster and a number of outliers on one or both sides of the cluster. Here, the overall condition number of the matrix is typically a very poor indicator of the complexity of solving the corresponding linear system, and iterative methods can still be very effective even if that condition number is very large.

\begin{definition}
Consider an $n\times n$ matrix $A$ with singular values $\sigma_1\geq \sigma_2\geq ...\geq \sigma_n>0$, and let $k,l\in\{0,...,n\}$ satisfy $k+l<n$. We say that $A$ is a $(n,k,l,\kappa)$ matrix if $\sigma_{k+1}/\sigma_{n-l}\leq \kappa$. 
\end{definition}
In other words, we split the singular values of $A$ into three groups: the top $k$, the bottom $l$, and the cluster with condition number $\kappa$ (other notions of condition number can be considered for the cluster, see below). We can choose either $k$ or $l$ to be 0 in order to focus on one-sided outliers. Early example of considering this model for solving linear systems with Krylov methods is \cite{axelsson1986rate}, where they study convergence of Conjugate Gradient on positive definite systems (equivalently, convergence of LSQR for general systems) in exact arithmetic.

\emph{Large outlying singular values.} 
We first focus on the setting of large outlying singular values (i.e., moderately large $k$, and $l=0$). This setting is typically motivated by problems arising in data science or statistics, where the input data comes from a low-dimensional manifold, which gives the matrix a low-rank structure that results in large  outlying singular values. The remaining cluster can be the result of noise in the data, or deliberate regularization of the system, i.e., $(A+\lambda I)x=b$ or $(A^\top A+\lambda I)x = b$ \cite{derezinski2024solving,derezinski2025fine,derezinski2024faster}.

\emph{Direct vs matrix-vector product complexity.} The large outlying singular values setting leads to a complexity separation between two standard models of computation.
\begin{enumerate}
    \item \textit{Direct access model:} We can perform computations with $A$ directly.
    \item \textit{Matrix-vector product model:} We interact with $A$ only via queries $v\rightarrow \{Av,A^\top v\}$.
    \end{enumerate}
In the direct access model, \cite{derezinski2025approaching} showed that we can solve an $(n,k,0,\kappa)$ linear system to high precision via randomized sketching techniques using
\begin{align*}
    O(k^\omega) + \tilde O(n^2\kappa)
\end{align*}
arithmetic operations (in exact arithmetic), where $\omega$ denotes the exponent of a matrix multiplication algorithm of your choice. On the other hand, in the matrix-vector product model, Krylov subspace methods need $\tilde O(k + \kappa)$ matrix-vector product queries to solve this system, and \cite{derezinski2025fine} showed that $\Omega(k)$ queries is necessary in this model. If the linear system is dense, this  leads to an additional cost of $\Omega(n^2k)$ operations, leading to a complexity separation for $k\gg \kappa$ (even if we restrict to classical matrix multiplication with $\omega=3$). 

\emph{Complexity for sparse matrices.} Importantly, the matrix-vector product model is most relevant for sparse matrices, in which case the extra cost paid by the Krylov solver is only $\tilde O(\mathrm{nnz}(A)\cdot k)$, where $\mathrm{nnz}(A)$ is the number of non-zero entries in $A$ (although, importantly, additional costs arise in finite precision arithmetic). Many randomized linear algebra methods can also leverage input sparsity for faster running time, so it is natural to expect that a similar improvement may be attainable in the direct access model.
\begin{problem}
    Can we solve an $(n,k,0,\kappa)$ linear system $Ax=b$ to high accuracy using $\tilde O(k^\omega + \mathrm{nnz}(A)\kappa)$ arithmetic operations?
\end{problem}

\emph{Averaged condition numbers.} Another way of exploiting large outlying singular values or eigenvalues is through averaged condition numbers. For illustration, consider replacing the usual condition number $\kappa = \frac{\sigma_1}{\sigma_n}$ with $\bar\kappa=\frac1n\sum_i\frac{\sigma_i}{\sigma_n}$. Clearly $\bar\kappa\leq \kappa$, and the gap can be significant, e.g., when $\sigma_i= 1/i$, then $\kappa= n$ but $\bar\kappa=O(\log n)$. In an early example of this, \cite{lee2013efficient} showed that we can replace $\kappa$ with $\bar\kappa$ in the complexity of solving psd linear systems when using randomized row-sampling methods rather than Krylov methods. This leads to the following more general definition \cite{derezinski2025approaching} (among others, e.g., see \cite{agarwal2020leverage}).
\begin{definition}\label{d:kappabar}
For $p \in (0, \infty]$, we define the $p$-averaged condition number of a sequence $\sigma_1\geq...\geq\sigma_n>0$  as
\begin{align*}
    \bar\kappa_{p} := \bigg(\frac1{n}\sum_{i}\frac{\sigma_i^p}{\sigma_n^p}\bigg)^{1/p}
    \text{ if } p \neq \infty
    \text{ and }
    \bar\kappa_{\infty} := \lim_{p \rightarrow \infty}  
    \bar\kappa_{p} 
    = \frac{\sigma_{1}}{\sigma_n}
    \,.
\end{align*}
\end{definition}
Naturally, one can incorporate the notion of an averaged condition number into the definition of $(n,k,l,\kappa)$ linear systems by excluding the top/bottom singular values from Definition \ref{d:kappabar}, but we will not do that for the sake of simplicity. Note that $\bar\kappa_{p_1}\leq \bar\kappa_{p_2}$ whenever $p_1<p_2$, which motivates the following question.
\begin{problem}\label{p:kappabar}
    Find the smallest $p>0$ such that we can solve an $n\times n$ linear system $Ax=b$ to high accuracy using $\tilde O(\mathcal{T}_A\cdot\bar\kappa_p)$ arithmetic operations, where $\bar\kappa_p$ is the $p$-averaged condition number of the singular values of $A$, whereas $\mathcal{T}_A$ is either $n^2$ or $\mathrm{nnz}(A)$.
\end{problem}
For $\mathcal{T}_A=n^2$, \cite{derezinski2025approaching} obtained this for any $p>1/2$ for general linear systems, and also improved to $p>1/4$ for psd linear systems. For $\mathcal{T}_A=\mathrm{nnz}(A)$, it is unclear that any such result is possible for $p<\infty$. However, partial results are known for $\mathcal{T}_A = n\cdot \max_i \mathrm{nnz}(A(i,:))$, which is proportional to $\mathrm{nnz}(A)$ when $A$ is uniformly sparse, e.g., \cite{lee2013efficient}. Problem \ref{p:kappabar} is closely connected with the complexity of Schatten $p$-norm estimation \cite{musco2018spectrum,derezinski2025approaching}.

\emph{Small outlying singular values.} We now turn to the more challenging setting of small outlying singular values (i.e., moderately large $l$, and $k=0$). This setting can be motivated by applications in PDEs where the eigenvalues often cluster at the top of the spectrum, and then decay towards the end. Small outlying singular values also naturally occur in certain random matrices, as discussed below. \cite{axelsson1986rate} showed that Krylov subspace methods can in exact arithmetic solve a $(n,0,l,\kappa)$ linear system to high accuracy using $\tilde O(l\kappa)$ matrix vector products, compared to $\tilde O(k+\kappa)$ for $(n,k,0,\kappa)$ systems. This suggests that capturing small outlying singular values by linear solvers is harder than capturing large ones. Moreover, there appears to be no existing work improving on this in the direct access model, for either dense or sparse matrices. This raises the following question.
\begin{problem}
    What is the complexity of solving $(n,0,l,\kappa)$ linear systems, and how is this affected by the model of computation (direct access vs matrix-vector product)?
\end{problem}

\emph{Solving a random linear system.}\footnote{Posed by Santosh Vempala at the workshop.} A natural test-case for small outlying singular values is a linear system based on a random $n\times n$ matrix. Consider for example a matrix with i.i.d. Gaussian entries. Its singular values decay linearly from $\sim\!\sqrt n$ to $\sim\!\frac1{\sqrt n}$, and so, the condition number of this matrix is $O(n)$, but it improves to $O(n/l)$ if we exclude the bottom $l$ singular values. Similar statements hold for other random matrices, e.g., with independent Bernoulli or Rademacher entries. Thus, improving the complexity of solving $(n,0,l,\kappa)$ systems is closely related to the following question.
\begin{problem}
    Let $A$ be an $n\times n$ random matrix with i.i.d. Gaussian (or Rademacher, or Bernoulli) entries, and take any $n$-dimensional vector $b$. Is it possible (conditioned on a high-probability event) to solve the linear system $Ax=b$ to high accuracy faster than $O(n^\omega)$?
\end{problem}

Note that this question is open (informally speaking) even if we restrict to classical matrix multiplication and let $\omega=3$. A related (and potentially easier) question is whether the matrix-vector product complexity for this problem is $\Omega(n)$.

\emph{Stability and finite-precision arithmetic.} All of the above statements and questions are posed under the Real-RAM model, where one performs all operations in exact arithmetic (infinite precision). However, the Krylov guarantees for outlying singular values (such as those of \cite{axelsson1986rate}) are prone to instability unless one uses extremely high precision, which significantly affects their bit complexity in the Word-RAM model \cite{greenbaum_89,MMS2018}. Similar concerns have been raised with regards to randomized sketching algorithms, particularly in the context of their backward stability (which tends to be a more restrictive notion of stability) \cite{EMN25}. In each case, these issues have been (at least partially) addressed in practice, but are not well understood theoretically. This motivates a rigorous analysis of bit complexity and backward stability of solving $(n,k,l,\kappa)$ linear systems.
\begin{problem} \label{problem29}
    What is the bit complexity of solving $(n,k,l,\kappa)$ linear systems? In particular, does the $\tilde O(k^\omega + n^2\kappa)$ running time for $(n,k,0,\kappa)$ systems extend to the Word-RAM model?
\end{problem}
\begin{update}
The reference
\begin{quote}
        Yang P. Liu, Hoai-An Nguyen, Richard Peng, Junzhao Yang, Faster Solvers for Sparse Systems with Large Spectral Outliers, 2026.
\url{https://yangpliu.github.io/pdf/faster-solvers-sparse-spectral-outliers.pdf}
\end{quote}
claims progress on Problem~\ref{problem29}.
For dense square normal equations $B^\top B x=b$ with explicit access to $B$, and
$\sigma_{k+1}(B)/\sigma_n(B)=O(1)$, it gives a 
randomized algorithm with bit complexity
$
n^{o(1)}\widetilde O\left(k^{\omega+\eta}+n^2\right)
$
for every fixed $\eta > 0$. Thus, the target complexity in
Problem~\ref{problem29} for $\kappa=O(1)$ is nearly attained.
The more general case, where only the linear system matrix $A$ is given and/or $\kappa$ need not be constant, remains open.
\end{update}
\medskip 

An important case of Problem~\ref{problem29}, which highlights the instability issues in Krylov subspace methods, is when we treat all of the singular values as ``outlying''. In this case, instead of solving the system directly in $O(n^\omega)$ time, one can choose to run $n$ iterations of a Krylov method such as the Conjugate Gradient / Lanczos algorithm. If the matrix $A$ is sparse, then in Real-RAM model this gives $\tilde O(n\cdot \text{nnz}(A))$ running time, potentially better than a direct solve. However, this is not true in the Word-RAM model, since the required precision grows with the rank of the Krylov subspace. This raises the following question.
\begin{problem}
    What is the bit complexity of solving an $n\times n$ linear system with $\text{nnz}(A)$ non-zeros? In particular, does the $\tilde O(n\cdot \text{nnz}(A))$ runtime extend to the Word-RAM model.
\end{problem}
\cite{peng2021solving} made progress towards this problem by analyzing the bit complexity of a block Krylov algorithm. Their argument, which was later improved by \cite{nie2022matrix}, showed that one can attain a running time faster than $O(n^\omega)$ for a sparse linear system in the Word-RAM model. However, the resulting complexity, $\tilde O(\text{nnz}(A)^{\frac{\omega-2}{\omega-1}}n^2)$, is still worse  than what is attainable in the Real-RAM model.

\subsection{Preconditioning and convergence of iterative solvers}

\scribe{Chen Greif, Jörg Liesen, and Richard Peng} 

The goal of preconditioning is to speed up the convergence of iterative methods such as CG and GMRES for solving linear algebraic systems $Ax=b$, particularly for large and difficult systems such as those arising regularly as discretizations of partial differential equations (PDEs). This is achieved by replacing the given linear system $Ax=b$ by a system of the form $$ M^{-1} A x = M^{-1} b,$$
for example, in the hope that the convergence of the solver is faster for this modified (preconditioned) system. Heuristically, $M$ should be a good approximation of $A$ (or rather, its inverse $M^{-1}$ a good approximation of $A^{-1}$), and at the same time the action of $M^{-1}$ must be easy to compute.

The above form of preconditioning is called {\em left preconditioning}, because we multiply the system by $M^{-1}$ on the left. Other forms are {\em right preconditioning}, where we form $$ A M^{-1} \tilde {x} = b, \qquad \tilde{x}= M x, $$ or {\em split preconditioning}, where we form 
$$ M_1^{-1} A M_2^{-1} \tilde{x} = M_1^{-1} b, \qquad \tilde{x}= M_2 x.$$
The differences between the effects of each preconditioning approach are interesting, but are beyond the scope of this document.

In general, finding an effective preconditioner is highly problem-dependent and is sometimes viewed as ``a combination of art and science''; see, e.g., \cite[Chapter~10]{Saad2003}. The well-known convergence bound of the CG method for SPD systems $Ax=b$, given by
$$\frac{\|x-x_k\|_A}{\|x-x_0\|_A}\leq 2\left(\frac{\sqrt{\kappa(A)}-1}{\sqrt{\kappa(A)}+1}\right)^k,\quad \kappa(A)=\frac{\lambda_{\max}(A)}{\lambda_{\min}(A)},$$
seems to indicate that a smaller condition number of the system matrix leads to faster convergence of the method. Therefore, it may seem that preconditioning techniques should aim mainly at reducing the condition number of the system matrix, i.e., achieving $\kappa(M^{-1}A)<\kappa(A)$. However, the above is just an upper bound. During the iteration, the CG method adapts in a nonlinear way to the location of the individual eigenvalues, i.e., its behavior depends on the entire spectrum of $A$. In general, a smaller condition number does not necessarily imply a faster convergence of the CG method! It is therefore not recommended to use the reduction of the condition number as the only criterion for the construction of preconditioners; see, e.g.,~\cite[Section~2.4]{CarLieStr24} for examples and further discussion of this essential point.

Another reason for focusing on condition number reduction rather than ``improving'' the entire eigenvalue distribution by preconditioning is that the eigenvalue distribution is usually unknown in practical applications. Yet, the recent line of work in~\cite{GerMadNieStr19,GerNieStr20,GerNieStr22,NieStr24} has shown that for certain elliptic PDE model problems, for example $-\nabla\cdot (k(x)\nabla u)=f$ with  $k(x)$ being  a measurable, bounded, and uniformly positive function from ${\mathbb R}^d$ to ${\mathbb R}$, one can obtain the exact localization of the eigenvalues of the discretized operators (i.e., the matrices), and of the left-preconditioned matrices under Laplacian preconditioning. This approach yields new insights into the convergence of CG for the preconditioned systems. 

\begin{problem}
Several specific open problems that exist in this context are mentioned in~\cite[Section~9]{NieStr24}.   More generally, try to extend the approach in~\cite{GerMadNieStr19,GerNieStr20,GerNieStr22,NieStr24} to other types of PDEs or preconditioners.   
\end{problem}

The behavior of iterative methods such as GMRES for solving nonsymmetric (more precisely, nonnormal) linear algebraic systems is much less understood than the behavior of CG. In general, the performance in this case may depend not only on the eigenvalue distribution, but also on further quantities including the conditioning of the eigenvectors of $A$, which motivates considering $\epsilon$-pseudospectra or the field of values of $A$ in the convergence analysis; see, e.g.,~\cite[Section~5.7]{LieStr13} for a survey of results. Despite numerous existing approaches, the following problem is largely open (cf. also~\cite{Emb22}):

\begin{problem}
    Develop a descriptive convergence bound for GMRES applied to nonnormal linear algebraic systems, that particularly describes the (transient) behavior early in the iteration. 
\end{problem}

Although the relationship between eigenvalues and convergence is in general unclear, many different preconditioning techniques in the context of nonsymmetric linear algebraic systems focus entirely on the clustering of eigenvalues. Authors then argue handwavingly that a few eigenvalue clusters mean that GMRES will converge in only a few steps. In order to make such arguments rigorous, at least in some cases, we propose the following problem statement.

\begin{problem}
    Find an interesting (nonselfdajoint resp. nonsymmetric) model problem, where the theoretical result on the eigenvalue distribution of the preconditioned matrix indeed corresponds to the actual convergence behavior of GMRES (e.g., $M^{-1}A$ has $m$ eigenvalue clusters and an acceptable approximate solution is obtained by GMRES in ${\cal O}(m)$ steps). Determine which particular properties of the problem implied this correspondence.    
\end{problem}

Researchers looking for a more specific open problem in this context may consider the following: The analysis of the GMRES convergence behavior for an SUPG discretized convection-diffusion model problem in~\cite{LieStr05} largely explains how the length of the initial phase of slow convergence depends on the boundary conditions of the PDE. It does not explain, however, which properties of the PDE, its discretization, or the linear algebraic system determine the convergence acceleration after the initial phase; see~\cite[Section~3.4]{LieStr05} for a more detailed discussion. We point out that the SUPG discretized operator is highly nonnormal in the sense that it has very ill conditioned eigenvectors; see, e.g.,~\cite[Section~4]{LieStr05} as well as~\cite{RedTre94} for a spectral analysis of a simple 1D conection-diffusion operator. The standard convergence bound on the GMRES relative residual norms, which contains the eigenvector condition number as a constant, cannot be descriptive in such cases; see, e.g., the discussion of this bound in~\cite[Section~3]{CarLieStr24}. One then has to resort to other  tools, such as field-of-values analysis \cite[Section 2.2]{Emb22}. We note, however, that while such an analysis is capable of establishing scalability in some relevant cases, it typically does not provide tight estimates of iteration counts.

\begin{problem}
Explain the convergence acceleration of GMRES after the initial (transient) phase that is observed in some convection-diffusion model problems such as the one in~\cite{LieStr05}. In particular, explain whether or not the speed of convergence can be related to eigeninformation (i.e., eigenvalue location and eigenvector conditioning) of the discretized convection-diffusion operators.    
\end{problem}

One of the most important goals in the context of linear systems that arise from discretizations of PDEs is to generate fully {\em scalable} preconditioned iterative solvers, i.e., solution procedures that produce nearly the same number of iterations regardless of the matrix dimensions and the parameters of the problem. In practical applications we are interested in a {\em small} number of iterations relative to the matrix size, and a {\em reasonable} error norm or residual norm (by reasonable we actually mean not too small, i.e., not close to machine unit roundoff). Many results for specific problems exist, but this is an open-ended topic and we therefore propose the following (cf. also Section~\ref{sec:scalable}):
\medskip
 
\textbf{Research direction} (Scalability with respect to physical parameters).
Design a fully robust  scalable (problem-dependent) preconditioner that is (almost) insensitive to changes in the values of physical parameters, e.g., small viscosity coefficients in fluid flow.

\subsection{Symmetric Krylov methods in finite precision arithmetic}

\scribe{Tyler Chen and Anne Greenbaum}

For symmetric matrices, Krylov subspace methods based on short recurrences (e.g. Lanczos or CG) are used to avoid orthogonalization costs. 
However, the behavior of such methods in finite precision arithmetic can be very different than in exact arithmetic; see \cite[Chapter 4]{chen_24} for a basic overview.
What can be said about the convergence rate and the ultimately attainable accuracy of these methods when used to solve symmetric positive definite linear systems or to evaluate a function of the matrix times a vector?
For some results on these topics, see, for example, \cite{druskin_knizhnerman_91,DGK98,greenbaum_89,Greenbaum97,MMS2018,paige_76,paige_80}.

A number of methods for solving linear systems (with or without a preconditioner) can be implemented in the following form:  Starting with an initial guess $x_0$ for the solution to $Ax=b$, for $k=1,2, \ldots$,
set 
\begin{eqnarray*}
x_k & = & x_{k-1} + a_{k-1} p_{k-1} \\
r_k & = & r_{k-1} - a_{k-1} A p_{k-1} .
\end{eqnarray*}
The updated vectors $r_k$ are not the true residuals $b-A x_k$ when these formulas are implemented in finite precision arithmetic.
It is pointed out in \cite{Greenbaum97} that these vectors often shrink well below the machine precision, and this is used to prove results about the size of the true residuals.
But this is not proved.

\begin{problem}
Taking the above recurrences to represent the CG or steepest descent algorithm, determine conditions on $A$ that will insure that the norms of the vectors $r_k$ drop below the machine precision in finite precision arithmetic or show some examples where they do not.
\end{problem}

The CG and Lanczos algorithms are equivalent in exact arithmetic, but they are implemented very differently.
Do they behave similarly in finite precision arithmetic?
In exact arithmetic, CG can be thought of as computing the Cholesky factorization of the tridiagonal matrix produced by the Lanczos algorithm.
But the relation between the two algorithms in finite precision arithmetic is still unexplained.

\begin{problem}
Show that the CG and Lanczos algorithms do or do not behave similarly in finite precision arithmetic when used to solve a symmetric positive definite linear system.  This will, of course, depend on the implementation of each method, so use standard good implementations of both.    
\end{problem}

A related problem is the following:

\begin{problem}
How many bits of precision are necessary to guarantee that CG applied to an $n$ by $n$ SPD linear system $Ax=b$ obtains an approximate solution with normwise backward error at most $\epsilon$ in $n$ (or fewer) steps?  
\end{problem}

Paige's analysis of the standard Lanczos algorithm in finite precision arithmetic  \cite[etc.]{paige_76,paige_80} describes the behavior of a standard implementation of the Lanczos algorithm in finite precision arithmetic. 
This analysis serves as the starting point for other analyses, such as Greenabum's analysis of CG \cite{greenbaum_89} and Knizhnerman's analysis of the Chebyshev moments \cite{knizhnerman_96}, Druskin and Knizhnerman's anlaysis of matrix function computations \cite{druskin_knizhnerman_91} as well as the recent analysis of Musco, Musco, and Sidford \cite{MMS2018}, and Druskin, Greenbaum, and Knizhnerman's earlier results on convergence rates and attainable accuracy of the method for solving symmetric positive definite linear systems \cite{DGK98}.

Paige's analysis of Lanczos is based on the standard model of floating point arithmetic, and does not analyze the impact of a preconditioner, which may be a matrix stored in memory or may be represented as an implicit operator that is applied inexactly.
Moreover, by switching to a more general model of matrix-vector multiplication and application of a preconditioner (e.g., $\mbox{fl}(Ax) = (A + \Delta A)x$, fl(solution to $Mz=r$) satisfies $(M + \Delta M)z = r$, without explicitly stating what $\Delta A$ and $\Delta M$ are), it may be possible to obtain a more intuitive/interpretable analog of existing analyses for finite precision arithmetic.

\begin{problem}
Analyze the Lanczos algorithm under general models of inexact matrix-vector multiplication and/or the use of a preconditioner.
Progress is made if these bounds generalize Paige's analysis and/or provide new, easier to understand, insight into the behavior of the algorithm.
\end{problem} 

In particular, one might consider an analysis of the infinite dimensional analog of Lanczos for computing orthogonal polynomials (often called the Stieltjes procedure).
By more explicitly leveraging the intrinsic connection between the Lanczos algorithm and orthogonal polynomials, 
one may gain further insight into the behavior of Lanczos-like methods.
For instance, this approach reveals that, in a certain sense, the forward instability of the Lanczos algorithm is due to the ill-conditioning of the bijection between moments of a certain distribution depending on $(A,b)$ and the tridiagonal Lanczos matrix \cite{knizhnerman_96,chen_trogdon_24,chen_24}.

Another natural generalization of Paige's analysis is to block Krylov methods.
So far, no ``Paige style'' analysis of the block Lanczos algorithm has been performed. 
A major challenge for analyzing block Krylov methods is that, even in exact arithmetic, the behavior can be very complicated. 
For example, the blocks may become rank deficient, and need to be deflated. 
Nevertheless, one hopes for an analysis of the block Lanczos algorithm that can provide theoretical justification and insight into the use of many block Krylov methods in finite precision arithmetic. 

\begin{problem}
    Extend Paige's analysis of the Lanczos algorithm \cite{paige_76,paige_80} to the block Lanczos algorithm.
    The exact analysis will depend on the use of an orthogonalization scheme used within blocks, so either specify certain assumptions on the accuracy of the scheme, or analyze a standard ``good'' algorithm.
\end{problem}

\subsection{The Forsythe conjecture}

\scribe{J\"org Liesen}

\emph{Background.} In a paper from 1968~\cite{For68}, Forsythe studied the asymptotic behavior of the ``optimum $s$-gradient method''. In modern terms he analyzed the asymptotic behavior of the CG method when it is restarted every~$s$ steps. He proved for the case of restart length $s=1$ (when the method is mathematically equivalent to steepest descent), and observed numerically for $s=2$, that the normalized error vectors (or normalized residual vectors) eventually cycle back and forth between two limiting directions. Such oscillatory behavior is unwanted, since it asymptotically leads to slow (at best linear) convergence of the method. Numerous authors have studied this behavior in the context of gradient descent methods for minimizing quadratic functionals; see, e.g., the survey~\cite{ZouMag22}. Forsythe conjectured that oscillatory limiting behavior occurs for every fixed restart length $s$. Several proofs of this behavior for $s=1$ have been published in addition to Forsythe's, and the case $s=2$ was considered in~\cite{ZhuBon83}. In general, however, the Forsythe conjecture remains largely open, which is quite surprising in light of the popularity and widespread use of methods like steepest descent, gradient descent, and CG.

\medskip

\emph{Mathematical details.} Forsythe considered the ``optimum $s$-gradient method'' for minimizing
\begin{equation}\label{eqn:f}
f(x)=\frac12 x^TAx-x^Tb,
\end{equation}
where $A\in {\mathbb R}^{n\times n}$ is symmetric positive definite and $b\in {\mathbb R}^n$. The unique minimizer of~$f$ is given by $x=A^{-1}b$, i.e., by the uniquely determined solution of the linear algebraic system $Ax=b$. Using modern notation, the uniquely determined iterate $x_1$ of the ``optimum $s$-gradient method'' starting with $x_0$ is defined by
\begin{equation}\label{eqn:forsythe}
x_1\in x_0+\K_s(A,r_0)\quad\mbox{such that}\quad x-x_1\perp_A \K_s(A,r_0),
\end{equation}
where $r_0=b-Ax_0$ and $\K_s(A,r_0)={\rm span}\{r_0,\dots,A^{s-1}r_0\}$. Note that the orthogonality condition in \eqref{eqn:forsythe} can also be written as $r_1\perp \K_s(A,r_0)$. Equivalently, the iterate $x_1$ satisfies
$$\|x-x_1\|_A=\min_{z\in x_0+\K_s(A,r_0)}\|x-z\|_A.$$
Thus, in exact arithmetic, $x_1$ is the uniquely determined $s$th iterate of the CG method applied to $Ax=b$ with the initial vector $x_0$. In order to avoid trivialities we assume that $1\leq s<d(A)$, where $d(A)$ is the degree of the minimal polynomial of $A$, as well as $d(A,r_0)\geq s+1$, where $d(A,r_0)$ is the grade of $r_0$ with respect to $A$. The last condition guarantees that $r_1=b-Ax_1\neq 0$. In fact, it can be shown that $d(A,r_0)\geq s+1$ implies that also $d(A,r_1)\geq s+1$. We now consider $x_1$ as the new initial vector and apply \eqref{eqn:forsythe} again to obtain $x_2$, i.e.,
$$x_2\in x_1+\K_s(A,r_1)\quad\mbox{such that}\quad x-x_2\perp_A \K_s(A,r_1),$$
and so on. This restarted method will converge to $x=A^{-1}b$. 

The Forsythe conjecture is concerned with the \emph{limiting directions} of the residuals in this restarted iteration. In order to study this behavior we consider a sequence of vectors constructed by restarting \eqref{eqn:forsythe} with an additional normalization:
\begin{align}
& \mbox{For \ensuremath{k=0,1,2,\dots}}\nonumber\\
& \hspace{1cm}\mbox{$y_k=r_k/\|r_k\|$,\quad where $r_k=b-Ax_k$,}\label{eqn:yk}\\
& \hspace{1cm}\mbox{$x_{k+1}\in x_k+\K_s(A,y_k)$
such that $x-x_{k+1}\perp_A \K_s(A,y_k)$.}\label{eqn:forsythe1}
\end{align}
As mentioned above, Forsythe proved that for $s=1$ the iteration in the limit cycles back and forth between two limiting directions, and he formulated the following conjecture: 

\begin{problem}[Forsythe conjecture] \label{problem220}
    For $2\leq s<d(A)$, each of the two subsequences $\{y_{2k}\}$ and $\{y_{2k+1}\}$ in \eqref{eqn:yk}--\eqref{eqn:forsythe1} has a single limit vector.
\end{problem}

Forsythe's conjecture is only about the existence of limit vectors, and not about the speed of convergence. However, if the conjecture holds, then asymptotically the vectors with the even indices become arbitrarily close to being collinear, and the same happens for the vectors with the odd indices. Therefore the error (or residual) norms of the restarted iteration can converge to zero at best linearly.

Numerical evidence suggests the conjecture is true. It can also be shown, for example, that
$$\lim_{k\rightarrow \infty} \|y_{2k+2}-y_{2k}\|=0$$
(and similarly for the odd indices), but this is not sufficient for proving the conjecture. A survey of the state-of-the-art of results in the context of the Forsythe conjecture until 2023 can be found in~\cite{FabLieTic23}. That paper also contains a generalization of the conjecture to nonsymmetric matrices and a ``cross iteration'' which alternates with $A$ and $A^T$. The ``cross iteration'' shows a similar limiting behavior as the iteration \eqref{eqn:yk}--\eqref{eqn:forsythe1}, which for certain orthogonal matrices and restart length $s=1$ is shown in~\cite{FabLieTic23}.

\medskip
\begin{update}
References
\begin{quote}
Matthew J. Colbrook, George Stepaniants and Alex Townsend, A Proof of the Forsythe Conjecture for the Two-Step Restarted Conjugate Gradient Method,
arXiv:2608.02852, 2026
\end{quote}

\begin{quote}
Jarek Liesen, Proof of the Forsythe Conjecture for $s = 2$, 2026.  \url{https://jarek.ai/papers/proof-of-the-forsythe-conjecture-for-s-2.pdf}.
\end{quote}
\begin{quote}
Richard Peng, Convergence of Normalized Residuals for Restarted Conjugate
Gradients of Length Two, 2026. \url{https://yangpliu.github.io/repo/restarted-cg-two/paper.pdf}
\end{quote}
claim to have solved Problem~\ref{problem220} for $s = 2$.
\end{update}

\section{Eigenvalue problems}
\label{sec:eig}
\subsection{Deterministic pseudospectral shattering / Minami bound}
\scribe{Peter Buergisser and Nikhil Srivastava} 

\noindent\motivation 
All known provable algorithms for computing a $\delta$-backward error diagonalization 
or Schur form of an arbitrary $A\in\mathbb{C}^{n\times n}$ begin by perturbing $A$ by a random matrix $E$ of size roughly $\delta$, and then running some algorithm on $A+E$. 
The purpose of this perturbation is to ensure that the condition number 
$\kappa_{eig}(A)$ of the diagonalization problem
becomes polynomially bounded in $n/\delta$ with high probability, i.e.,
\begin{equation}\label{eqn:shatteringbound}
\kappa_{eig}(A+E)=O((n/\delta)^c) ,
\end{equation} 
where $c$ denotes some constant.
This means that the distance of $A$ to the set of matrices with multiple eigenvalues is at least inverse polynomial in $n/\delta$.
A geometric definition of the condition number of the eigenpair can be found in 
\cite[\S14.3]{BurgisserCucker_Condition}. 
A direct definition is 
$$
 \kappa_{eig}(A) = \kappa_V(A)/gap(A) ,
$$
where $gap(A)$ denotes the minimum gap between two eigenvalues of $A$, and 
$\kappa_V(A)$ denotes the condition number of the matrix of normalized 
eigenvectors of $A$, which is well defined if $A$ has no multiple eigenvalues;
see~\cite{Banks-et-al-23}. 
The process of replacing $A$ by $A+E$ is called ``pseudospectral shattering'' due to its manifestation in the pseudospectrum.

As a result, all known provable algorithms for diagonalization/Schur form are randomized. As noted in \cite{sobczyk2024deterministic}, this motivates the question of finding deterministic algorithms, and a natural path towards such algorithms is to derandomize the regularizing perturbation $E$.

\begin{problem}[Deterministic Pseudospectral Shattering]
Give a provably correct and efficient \emph{deterministic} algorithm for the following task: on input $A\in\mathbb{C}^{n\times n}$ and $\delta>0$, output a perturbation $E$ with $\|E\|\le \delta$ such that $\kappa_{eig}(A+E)\le C(\delta/n)^c$,
where $C,c$ are universal constants. 
In this context, efficient means $O(n^3\log^c(n/\delta))$ arithmetic operations, either in exact arithmetic or, more ambitiously, in finite arithmetic with $O(\log^c(n/\delta))$ bits of precision.
\end{problem}

The paper~\cite{Banks-et-al-23} describes a randomized algorithm achieving such 
pseudospectral shattering. The point of the above problem is whether 
a \emph{deterministic} algorithm exists, satisfying the same requirements on 
time and precision.

\remarks
\begin{itemize}
    \item Davies \cite{davies2008approximate} shows how to deterministically obtain a perturbation $E$ that yields $\kappa_V(A+E)\le (n/\delta)^{n}$. 
    \item There are two reasons for wanting a polynomial bound on $\kappa_{eig}$, 
    as opposed to an exponential one. First, the iteration count of algorithms depends polynomially on $\kappa_{eig}$ (as in \cite{Armentano-et-al}) or logarithmically 
    on $\kappa_{eig}$ as in \cite{Banks-et-al-23}. Second, the number of bits of precision required by the algorithms depends logarithmically 
    (as in \cite{Armentano-et-al}) or polylogarithmically (as in \cite{Banks-et-al-23}) on $\kappa_{eig}$. This is not just an artifact of the proofs, as it is easy to see that the power method can  converge slowly on matrices 
    with large $\kappa_{eig}$, suggesting slow convergence of other algorithms as well.
    
    \item For random dense complex Gaussian $E$, it is known \cite{Banks-et-al-23} that one obtains \eqref{eqn:shatteringbound} with $c\le 5$ with high probability (the reason for not specifying $c$ precisely is that it depends on how high a probability one wants).
    \item For sparse complex Gaussian $E$ where each entry is nonzero with probability $O(\log^2(n)/n)$, the weaker bound of $\kappa_{eig}(A+E)\le \exp(\log^2(n/\delta))$ is known to hold with high probability \cite{sparsepseudospectralshattering}, which is worse but still useful in algorithms which depend logarithmically on $\kappa_{eig}$. 
    \item (On derandomization) 
    Homotopy methods for solving systems of polynomial equations require well-conditioned starting systems. Is is known that random systems are well-conditioned, but constructing them efficiently is difficult. 
\cite{etayo-et-al} solved this problem for univariate polynomials. 
Generally, derandomization is a major topic of study in theoretical computer science. 

\end{itemize}
There is a related problem is about diagonal random perturbations of Hermitian matrices, which are useful in Hermitian diagonalization \cite{Banks-et-al-23, shah2025hermitian} and in Krylov subspace methods \cite{meyer2024unreasonable, kressner2024randomized}.

\begin{problem}[Deterministic Minami Perturbation]
    Minami \cite{minami1996local} showed that if $A$ is $n\times n$ Hermitian tridiagonal and $E$ is a random diagonal matrix with independent absolutely continuous entries, then $gap(A+\delta E)\ge C(\delta/n)^c$ for some constants $C,c$ depending on the distribution of the entries. Give an efficient deterministic algorithm for finding such a perturbation. Here, efficient means faster than the parent algorithm in which it is used, so $O(n^3\log^c(n/\delta))$ for diagonalization and, ambitiously, $O(nnz(A)\log^c(n/\delta))$ or $O(nnz(A)\log^c(n)/\delta^c)$ for Krylov methods.
\end{problem}

A possible approach is to start with tridiagonal Toep\-litz matrices, for which an explicit formula is known for the eigenvalues.

\subsection{Accuracy vs. precision in efficient diagonalization algorithms}
\scribe{Peter Buergisser, Kate Pearce, and Ryan Schneider}

We propose the following problem, which will improve the result in~\cite{Banks-et-al-23} (for diagonalization) and in \cite{banks2022global} (for the Schur form). The motivation is that the amount of precision required by the best currently known provable algorithms for diagonalization/Schur form is much larger than the gold standard of $\log(1/\delta)+c\log(n)$ bits preferred in numerical analysis (which has been achieved for the Hermitian diagonalization problem \cite{shah2025hermitian}).

\begin{problem}
    Give a provable randomized algorithm for backward error diagonalization or Schur form with runtime $O(n^3\log^c(n/\delta))$ that only uses $O(\log(n/\delta))$ bits of precision.
\end{problem}

By contrast, the corresponding result in~\cite{Banks-et-al-23} runs in nearly matrix multiplication time, but uses $O(\log^4(n/\delta)\log(n))$ bits of precision, and the algorithm of \cite{banks2022global} runs in $O^*(n^3\log^2(n))$ operations with $O^*(\log^2(n/\delta))$ bits of precision plus $O^*(n\log^c(n))$ operations with $O^*(\log^4(n/\delta))$ bits of precision, where the $O^*$ notation suppresses $\log\log$ factors. The algorithm in \cite{Armentano-et-al} only uses $O(\log(n/\delta))$ bits of precision, although the running time is much worse (but still polynomial in $n$).

\remarks 
    The additional bits of precision required in \cite{Banks-et-al-23} stem from the use of inversion in the subroutine ``SGN'', which computes the matrix sign function via a standard Newton iteration. If this routine is replaced by an inverse-free method (e.g., the Implicit Repeated Squaring algorithm of Malyshev \cite{Malyshev1,Malyshev2}), then only the desired $O(\log(n/\delta))$ bits of precision are required (see the floating-point bounds in \cite[Chapter 6]{Schneider_Thesis} or \cite[Appendix B]{Demmel_deflating_subspace}. While a full end-to-end floating-point analysis remains open, this would imply an optimal algorithm for diagonalizing not only any matrix but any definite matrix pencil via \cite{Demmel2024}.

\subsection{Ritz values and approximate invariant subspaces}

\scribe{Rikhav Shah and Mark Embree}

We propose several problems related to the approximation of invariant subspaces and the eigenvalue approximations drawn from the Rayleigh--Ritz process using such subspaces.

\begin{problem} \label{prob:ritz1}
For a given $A\in\C^{n\times n}$, let $\widehat{\cal V}$ denote an approximate invariant subspace of dimension $d$ of $A$, in the sense that there exists some $d$-dimensional invariant subspace ${\cal V}$ of $A$ such that $\angle({\cal V},\widehat{\cal V})$ is small.  
Let ${\cal W}$ be some other subspace of dimension $r<n-d$.  
Suppose the columns of $Q$ form an orthonormal basis for the (generically) $d+r$ dimensional space 
${\cal Q} = \widehat{\cal V} + {\cal W}$.  

Provide conditions that will ensure that $d$ eigenvalues of $Q^*AQ$ approximate the $d$ eigenvalues of $A$ associated with the invariant subspace ${\cal V}$. Impose assumptions about the eigenvalues of $A$ corresponding to $\cal V$, as needed.
\end{problem}

\motivation\ \ The restarted Arnoldi process builds a Krylov subspace of dimension $d+r$, where $d$ is the desired number of eignvalues and $r$ is a buffer space used to design the polynomial filter for the restart process, whereby the Krylov starting vector $v_1$ is replaced by $v_1^+ = p(A)v_1$ for a polynomial whose roots (hopefully) are distributed across the unwanted part of the spectrum.  This restart process enriches the components of $v_1^+$ in the desired eigenvectors.  After numerous restart cycles, the Krylov subspace ${\cal K}_{d+r}(A,v_1^+)$ will (hopefully) contain a good approximation to a $d$-dimensional invariant subspace, extended by an additional $r$ dimensional space about which we can say nothing in particular.  In a sense, then, Problem~\ref{prob:ritz1} resembles  a local convergence theory for the restarted Arnoldi method.  The restarted Arnoldi convergence theory in~\cite{BER04} could potentially provide a framework.

The setting of Problem~\ref{prob:ritz1} could allow the implementation of heuristics with little-to-no drawback when the heuristics fail. In particular, if it is true that $Q^*AQ$ contains eigenvalues close to those corresponding to $\cal V$ irrespective of $\cal W$, then one is free to add any vectors to $\cal Q$ in the hopes that they help approximate an invariant subspace, without doing any harm.

Useful sources include Sorensen's work on the implicitly restarted Arnoldi method (which includes a convergence theorem for the Hermitian case)~\cite{Sor92}, Saad's text on large-scale eigenvalue algorithms~\cite{Saa11}, and papers on Arnoldi convergence theory such as~\cite{BSS10}.

\begin{problem} \label{prob:ritz2}
Let $Q$ be an orthonormal basis for a Krylov space for arbitrary $A\in\C^{n\times n}$ with random starting vector $b$. Is $\kappa_V(Q^*AQ)$ bounded (by a polynomial in $n$) with high probability? Or is $\E[\area\Lambda_\eps(Q^*AQ)]$ (the expected value of the area in $\C$ of the $\varepsilon$-pseudospectrum) bounded by $\poly(n)\eps^\beta$ for $\beta$ close to 2?
(This problem is open even when $A$ is a circulant shift matrix.)
\end{problem}

\motivation\ \ A straightforward calculation shows that if the range of $Q$ is close to containing an eigenvector of $A$, then the corresponding eigenvalue will be contained in $\Lambda_\eps(Q^*AQ)$. To conclude that this condition implies the Ritz value will approximate that eigenvalue well, one needs $\Lambda_\eps(Q^*AQ)$ to be small. Unfortunately, this need not be the case for all choices of $b$.  For example, if $A$ is the circulant shift matrix and $b = e_1$, then for all $k<n$, the Rayleigh quotient $Q^* A Q$ is a Jordan block (and thus $\Lambda_\eps(Q^*AQ)\propto\eps^{2/k}$ and $\kappa_V(Q^*AQ)=\infty$).  However, such examples are thought to be rare, and so this problem could potentially be avoided by a random selection of $b$.

Potentially useful resources include Stewart's backward error analysis for Krylov subspaces~\cite{Ste02}, the adversarial constructions of Duintjer Tebbens and Meurant~\cite{DM12}, and the analysis of $\E[\area\Lambda_\eps(Q^*AQ)]$ when $Q$ is Haar unitary in~\cite{Sha25}.

\medskip
\begin{update} Reference
\begin{quote}
Richard Peng. A Diagonalizable Obstruction for Random Two-Step Ritz Compressions, 2026. https://yangpliu.github.io/repo/two-step-ritz-obstruction/paper.pdf
\end{quote}
claims to provide negative answers to both questions of Problem~\ref{prob:ritz2} as stated. It constructs, for every even $n$, a real diagonalizable $A_n$, uniformly bounded in norm, such that for a standard real Gaussian starting vector $b$ and $k = 2$,
\[
\Pr\left\{
\kappa_V(Q^T A_n Q) \ge \frac{e^n}{5}
\right\}
\ge 1-2e^{-n/40}.
\]
Additionally, a proof is provided that, at $\varepsilon_n=5e^{-2n}$,
\[
\mathbb{E}\!\left[
\operatorname{area}\Lambda_{\varepsilon_n}(Q^T A_nQ)
\right]
\ge
\left(1-2e^{-n/40}\right)\frac{\pi e^{-2n}}{4},
\]
which rules out a uniform $\operatorname{poly}(n)\varepsilon^\beta$ estimate for every $\beta>1$. 

However, it should be noted that $\kappa_V(A_n)$ grows exponentially with $n$.
It remains open whether Problem~\ref{prob:ritz2} has a negative answer under the additional assumption that $A_n$ is normal or even under the weaker assumption  that $\kappa_V(A_n)$ grows at most polynomially with~$n$.
\end{update}

\begin{problem}  Consider a non-Hermitian $A\in\C^{n\times n}$ and suppose $Q\in\C^{n\times k}$ has orthonormal columns.  We know that the Ritz values (eigenvalues of $Q^*AQ$) are located in the numerical range, $W(A) = \{v^*Av: \|v\|_2=1\}$, a closed, convex subset of $\C$ that contains the eigenvalues of $A$ (and, potentially, points far from these eigenvalues).
\begin{itemize}
\item What can be said (deterministically or probabalistically) about the distribution of the Ritz values across $W(A)$ for $1<k<n$?
\item What if $Q$ is restricted to be a basis for a Krylov subspace?  What if $Q$ is the basis for a random subspace?
\end{itemize}
\end{problem}

\motivation\ \ The Ritz values are used as estimates for the eigenvalues of $A$.\ \ In the Hermitian case, Cauchy's interlacing theorem provides a complete answer to the deterministic case for arbitrary subspaces.  (Any arrangement of Ritz values that obeys interlacing is possible; see, e.g., \cite{PS08}.)  The non-Hermitian case seems highly challenging; the case $k=1$ is essentially solved (by definition, any $\theta \in W(A)$ has a generating unit vector $v$; the challenge is to find $v$, given $\theta$)~\cite{Car09b}, and $k=n$ (where the Ritz values must match the eigenvalues of $A$).  Results exist for the case of $k=n-1$ and normal $A$~\cite{CH13,Mal04}.  Some necessary conditions are provided in~\cite{CE12}.  Other sources include~\cite{DM12},\cite{GS12}, and~\cite{Sha25}.

\subsection{Revisiting \boldmath ${\rm MR}^3$}

\scribe{Ryan Schneider}

The Multiple Relatively Robust Representations (MRRR or MR$^3$) algorithm computes the eigenvalues and eigenvectors of a symmetric tridiagonal matrix $T$. Essentially a special version of bisection and inverse iteration, MR$^3$ can compute $k$ eigenpairs of $T$ in $O(kn)$ operations, and in this sense is optimal, i.e., it produces a full set of (numerically orthogonal) eigenvectors in $O(n^2)$ operations, which is a trivial lower bound. Accordingly, MR$^3$ yields the fastest algorithm for the symmetric eigenvalue problem via the following two-step procedure:
\begin{enumerate}
    \item Reduce the input (symmetric) matrix $A$ to a symmetric tridiagonal $T$ (e.g., by applying standard Householder reflectors).
    \item Diagonalize $T$ by applying MR$^3$.
\end{enumerate}

\indent MR$^3$ was introduced in \cite{Dhillon:CSD-97-971} and subsequently presented in two papers of Dhillon and Parlett \cite{Dhillon_orthogonal,DHILLON20041}. It was added to LAPACK shortly after \cite{Dhillon:CSD-04-1346}, and is now the default for the symmetric eigenvalue problem in Julia. 

There is currently no rigorous proof (in the vein of theoretical computer science) of success for a floating-point implementation of MR$^3$ when applied to an arbitrary symmetric tridiagonal matrix. Willems and Lang \cite{Willems_framework} have given a proof of correctness which is conditional on several underlying assumptions, some of which cannot be verified until the algorithm is executed on a given input. At the same time, it is known that MR$^3$ can fail in the presence of roundoff when eigenvalues are too tightly clustered, if implemented without care (see \cite{Parlett:CSD-04-1367}). 

\begin{problem}
    Give a ``TCS-style'' analysis of MR$^3$ which proves its correctness and convergence under easily verifiable conditions on the input.
\end{problem}

The second open problem presented in this section concerns the (as yet unsuccessful) application of this algorithm to the bidiagonal SVD problem.

\begin{problem}
    Let $B$ be an $n \times n$ (upper) bidiagonal matrix. Is there an algorithm that, in $O(kn)$ operations, can produce $k$ triples $(\sigma_i,v_i,u_i)$ such that $\sigma_i \geq 0$, $v_i, u_i \in {\mathbb C}^n$ are unit vectors, and, for machine precision $\epsilon$,
    \begin{enumerate}
        \item $\|Bv_i-\sigma_iu_i\|_2 = O(\epsilon n\|B\|_2)$;
        \item $|v_i^Hv_j|$, $|u_i^Hu_j| = O(\epsilon n)$ if $i \neq j$.
    \end{enumerate}
\end{problem}

In theory, MR$^3$ (which satisfies analogous error bounds for the symmetric tridiagonal eigenvalue problem) should answer this question in the affirmative. Naively, we might hope that by transforming the SVD problem $B = U\Sigma V^H$ into a companion tridiagonal eigenvalue problem (e.g., $B^HB = V\Sigma^2 V^H$, $BB^H = U\Sigma^2U^H$ or by looking at the $2n \times 2n$ Golub--Kahan matrix\footnote{The Golub--Kahan matrix is a permuted version of the $2n \times 2n$ block matrix  $\begin{pmatrix} 0 & B \\ B^T & 0 \end{pmatrix}$, which has eigenvector matrix $ \frac{1}{\sqrt{2}} \begin{pmatrix} U & U \\ -V & V \end{pmatrix}$.}) we could simply call MR$^3$ as a black box and be done. Nevertheless, this approach does not seem to work. If the left and right singular vectors $u_i$ and $v_i$ are computed independently by applying MR$^3$ to $BB^H$ and $B^HB$, the residual error $\|Bv_i-\sigma_iu_i\|_2$ can be large. On the other hand, working with the Golub--Kahan matrix often leads to insufficiently orthogonal singular vectors. \\
\indent These pitfalls were noticed by Großer and Lang \cite{GROER200345}, who proposed a modified version of MR$^3$ for the bidiagonal SVD that alleviates these issues (or at least appears to) by working with $B^HB$, $BB^H$, and the Golub--Kahan matrix simultaneously. While this algorithm was at one point proposed for LAPACK \cite{Willems_lapack}, it was eventually abandoned due to inaccuracies. The problem was taken up again about a decade later by Willems and Lang \cite{willems_etna,Willems_factorization,Willems_framework}, producing another implementation that was also ultimately ruled unreliable (see \cite[Appendix B]{Marques_SVD}). At the time of writing, this is the most recent work on MR$^3$. \\
\indent The question therefore remains:\ can we make a version of this algorithm work for the bidiagonal SVD?\ \ There are a number of potential directions to pursue here.
\begin{enumerate}
    \item \textbf{TCS Style Proof for the Symmetric Case.} Pursuing a full end-to-end proof for the real symmetric case may illuminate a gap in the existing theory that explains the failure in the bidiagonal setting (and could potentially be closed by, for example, randomizing). 
    \item \textbf{Understanding bidiagonal failure modes.} While proposed LAPACK implementations of MR$^3$ for the SVD problem have been rejected twice, little has been recorded about specific failure examples. A study aiming to extract the features of problematic bidaigonals (e.g., distributions of singular values) could motivate modifications to existing implementations. 
    \item \textbf{Debugging existing implementations.} Finally, it is unclear whether the implementation of Willems and Lang (tested in \cite{Marques_SVD}) fails because of a bug or because of a genuine gap in the theory.\footnote{This note is based on correspondence between Marques and Willems.} In-depth debugging is therefore warranted.
\end{enumerate}
\subsection{Largest / rightmost eigenvalue}
\scribe{Rikhav Shah}
\begin{problem}\label{problem:gapfreelambda}
Consider a general nonnormal matrix $A = V \Lambda V^{-1}$ with largest eigenvalue magnitude $\lambda_1$ and finite eigenvector condition number $\kappa_V = \|V\| \|V^{-1}\|$.
    Is there an algorithm that uses $\poly(\log(n\kappa_V),1/\epsilon)$ matrix-vector products with $A$ and returns $\lambda$ so that $|\lambda - \lambda'| \le \epsilon |\lambda_1|$, where $\lambda'$ is any eigenvalue of $A$ with $|\lambda' | \ge (1-\epsilon)\abs{\lambda_1}$?
\end{problem}

\remarks
It is possible to argue that power method gives an $\epsilon$ approximation to the spectral radius in a mixed forward-backward error sense using $O(\log(n\kappa_V)/\epsilon)$ iterations (see~\cite{sparsepseudospectralshattering}).
One candidate approach for upgrading this to a solution to the above problem is to ``search'' for an eigenvalue along the circle of radius the spectral radius found by the power method; namely, break the disk into $\eps$-sized chucks and run inverse iteration at a point in each chunk.

One strengthening of this problem is to demand an eigenvalue approximation for \textit{each} $\lambda'$ with $\abs{\lambda'}\ge(1-\epsilon)\abs{\lambda_1}$.

Finally, one can also ask a variant of the above question for approximately computing the rightmost eigenvalue of a matrix (in the complex plane). This was suggested by Lin Lin.

\section{Low-rank approximation and index selection}
\label{sec:lowrank}

In this section, we describe open directions for the column subset selection problem (CSSP)\footnote{Analogously, we may also be interested in selecting rows.} and computing rank-revealing factorizations. The CSSP is stated as follows: given a matrix $A \in \R^{m \times n}$, the goal is to select $1\leq k \leq \min(m, n)$ column indices $J = \{j_1, \dots, j_k\}$ such that $\Pi_J A$ is a good low-rank approximation to $A$, where $\Pi_J$ is the projection onto the column space of $A(:, J)$. There is no efficient algorithm that finds the optimal $J$, as this requires an exhaustive search through $\binom{n}{k}$ choices of columns. However, a number of algorithms can efficiently obtain a quasi-optimal choice of $J$, i.e., one satisfying 
\begin{equation}\label{eq:CSSP}
\|\Pi_J A - A \|_{\{2,F\}} \leq c(n, k) \min_{\text{rank}(B) \leq k} \|A - B\|_{\{2,F\}},
\end{equation}
where $c$ grows slowly (e.g., polynomially or logarithmically) with $n$ and $k$ and we may be interested in either the 2 or Frobenius norm.\footnote{Typically the natural $n$ dependence is a consequence of the norm chosen and the $k$ dependence is determined by the algorithm.} 

In a related problem of computing a rank revealing factorization, one wants to construct a rank-$k$ approximation $A_k$ by some algorithm. We call this algorithm a rank revealer if there exists a function $\mu(m, n, k)$, which is bounded by a polynomial in $m$, $n$, and $k$, such that~\cite{damle2025estimating}
\begin{equation}\label{eq:RRQR1}
    \frac{1}{\mu(m, n, k)} \sigma_i(A) \leq \sigma_i (A_k) \leq \mu(m, n, k) \sigma_i(A) \quad \text{for} \quad 1 \leq i \leq k
\end{equation}
and
\begin{equation}\label{eq:RRQR2}
    \frac{1}{\mu(m, n, k)}\sigma_{k + i}(A)\leq \sigma_{i}(A - A_k)  \leq \mu(m,n, k) \sigma_{k + i}(A) \quad \text{for} \quad 1 \leq i \leq \min(m, n) - k.
\end{equation}
The truncated singular value decomposition (SVD), obtained by retaining the
leading $k$ rank-one terms of the full SVD of $A$, provides the optimal
rank-$k$ approximation in the sense that it minimizes $\mu(m,n,k)$ and
achieves $\mu(m,n,k) = 1$. Nevertheless, the performance of a rank-revealing
method is not judged solely by this optimality criterion. In many applications,
additional considerations play a decisive role, including:
(1) the computational cost required to form $A_k$,
(2) the interpretability of the resulting column and row spaces, and
(3) the ability of the approximation to preserve structural properties of the
original matrix.

These considerations help explain the continued importance of
Gaussian elimination (GE) and QR-based rank-revealing techniques.
Compared with the SVD, such methods are often significantly cheaper to
compute, can produce low-rank approximations whose factors are easier to
interpret, and frequently retain useful structural features of the matrix.
Moreover, rank-revealing algorithms are commonly employed to construct
well-conditioned bases for the column space or row space of a matrix.

In practice, it is desirable to have $\mu(m,n,k)$ grow very slowly in $m$, $n$, and $k$. There are rank revealing algorithms in the literature that achieve $\mu(m,n,k)= \mathcal{O}(k\sqrt{mn})$ for GE-like factorization~\cite{miranian2003strong} and others that achieve $\mu(m,n,k)= \mathcal{O}(\sqrt{kn})$ for QR-like decompositions~\cite{gu1996efficient}. The rank-revealing QR factorization (RRQR) is one of the most popular rank-revealing algorithms~\cite{gu1996efficient}.

\subsection{When do greedy/randomized column selection algorithms work?} \scribe{David Persson, Diana Halikias, Paul Beckman, Zhen Huang, Anil Damle, Mark Fornace, and Alex Townsend}

Many practical algorithms for this problem incorporate randomization, adaptivity, and/or a greedy selection strategy. For example, the column-pivoted QR algorithm (CPQR)~\cite{businger1965linear}
is a greedy method that is observed to work very well on ``most'' matrices in the sense of~\ref{eq:RRQR1} and~\ref{eq:RRQR2}. However,  the worst case theoretical bound gives $\mu(n, k) = 2^k\sqrt{n-k},$ see, e.g., \cite{gu1996efficient}, indicating the potential for unfavorable exponential growth in 
$k$. How can we  understand this gap between theory and practice?

There is a similar story for  low-rank approximation using index selection. For LU with complete pivoting,  an exponential factor also appears in the optimality bound in general. But, the factor is only  polynomial if  the matrix is a  positive definite, Lipschitz-smooth kernel~\cite{jeong2025convergence}. One can also  reduce the exponential factor to polynomial in the maximum volume submatrix problem if the matrix is doubly diagonally-dominant~\cite{cortinovis2020maximum}.  In general, there is a sense that  general optimality bounds for these algorithms are quite pessimistic relative to how well they perform in practice. We seek to understand when and why these methods are provably quasi-optimal. 

\begin{problem}Under what structural or problem-dependent  assumptions on 
$A$ does there exist an efficient, greedy and/or randomized algorithm that  achieves quasi-optimal guarantees such as
\begin{enumerate}
\item the CSSP bound~\ref{eq:CSSP}, and
\item the RRQR bounds~\ref{eq:RRQR1}–\ref{eq:RRQR2}?
\end{enumerate}
Concretely, a solution to this problem would pick a structural assumption for $A$ and then provide an algorithm (or pick an existing one) that provably achieves near optimal bounds for either the CSSP or RRQR problem.
\end{problem}

\subsubsection{The discrete Lehmann representation}

\source{David Persson and Zhen Huang}

An important computational task in quantum physics is the approximation of the \emph{imaginary-time Green's function}
\begin{equation*}
    G(t) := \int_{-\Lambda}^{\Lambda} K(t,\omega) \rho(\omega) d\omega, \quad t \in [0,1],
\end{equation*}
where $K(t, \omega):= \frac{e^{-t\omega}}{1+e^{-\omega}}$ is the fermionic kernel, $\rho$ is an \emph{unknown} spectral density, and $\Lambda$ is a cutoff parameter. We are restricted to evaluating $G$ at a few nodes $t_1,\ldots,t_k$, which we are free to select, and we seek an approximation of the form
\begin{equation}\label{eq:dlr_approx}
    G(t) \approx \widehat{G}(t):=\sum\limits_{j=1}^k K(t,\omega_j) g_j, \quad G(t_i) = \widehat{G}(t_i). 
\end{equation}
for some suitably chosen $\omega_1,\ldots,\omega_k \in [-\Lambda,\Lambda]$ and coefficients $g_1,\ldots,g_k$ \cite{dlr}. 

The first condition in \eqref{eq:dlr_approx} states that we want to approximate $G$ with a linear combination of ``columns'' of the kernel $K$, while the second condition in \eqref{eq:dlr_approx} enforces interpolation at the chosen points $t_1,\ldots,t_k$. For fixed node sets $\{t_1,\ldots,t_k\}$ and $\{\omega_1,\ldots,\omega_k\}$, define the \emph{cross-approximation}\footnote{This assumes that the core matrix $\{K(t_i,\omega_j)\}_{i,j = 1}^k$ is invertible.}
\begin{equation*}
    \widehat{K}(t,\omega) := \begin{bmatrix} K(t,\omega_1) & \cdots & K(t,\omega_k) \end{bmatrix} \begin{bmatrix} 
    K(t_1 ,\omega_1) & \cdots & K(t_1,\omega_k)\\
    \vdots & \ddots & \vdots\\
    K(t_k ,\omega_1) & \cdots & K(t_k,\omega_k)
    \end{bmatrix}^{-1}  \begin{bmatrix} K(t_1,\omega) \\ \vdots \\ K(t_k,\omega) \end{bmatrix}.
\end{equation*}
One can verify that
\begin{equation*}
    \widehat{G}(t) = \int_{-\Lambda}^{\Lambda} \widehat{K}(t,\omega)\rho(\omega)d\omega,
\end{equation*}
and the approximation error satisfies
\begin{equation*}
    \|G-\widehat{G}\|_{\infty} \leq \|K-\widehat{K}\|_{\infty}\|\rho\|_1,
\end{equation*}
where $\|\cdot\|_p$ denotes the standard $L^p$-norm for functions. Hence, $\widehat{G}$ is an accurate approximation to $G$ if the nodes $\{t_1,\ldots,t_k\}$ and $\{\omega_1,\ldots,\omega_k\}$ form a good cross-approximation to $K$. 

A classical approach to selecting such nodes is Gaussian elimination with complete pivoting (GECP), a greedy method that attempts to maximize the volume of the cross-matrix $\{K(t_i,\omega_j)\}_{i,j = 1}^k$ \cite{bebendorf2000approximation,townsend2014computing}. Existing bounds \cite{cortinovis2020maximum,townsend2015continuous} rely on the assumption that the kernel is analytic in a sufficiently large Bernstein ellipse, and guarantee that GECP produces a cross-approximation to $K$ satisfying
\begin{equation}\label{eq:dlr_theorate}
    k = O\left(\Lambda + \log(1/\varepsilon)\right) \Rightarrow \|K-\widehat{K}\|_{\infty} \leq \varepsilon.
\end{equation}
However, one can prove the existence of a rank-$k$ approximation $K_k$ satisfying the sharper bound
\begin{equation*}
    k = O\left(\log(\Lambda)\log(1/\varepsilon)\right) \Rightarrow \|K-K_k\|_{\infty} \leq \varepsilon.
\end{equation*}
Empirically, the performance of GECP aligns much more closely with the latter rate than the theoretical rate given in \eqref{eq:dlr_theorate}. In applications, $\Lambda$ can be extremely large, say $\Lambda = 10^5$ or $10^6$, highlighting a substantial gap between the empirical performance of GECP and the existing theoretical guarantees. 

\begin{problem} \label{problem:gecp} Prove stronger theoretical bounds for GECP applied to $K$, by exploiting the structure of the kernel $K$. \end{problem}

Interesting, a relatively satisfaction answer to problem 4.2 is available when one assumes that $K$ is a continuous positive definite kernel. In this case, GECP simplifies to the pivoted Cholesky algorithm (also called the P-greedy algorithm in the kernel literature~\cite{santin6convergence}), where at each iteration the next step is the maximum pivot on the diagonal. Here, it is known that $\|K - K_k\|_\infty$ converges at a rate dictated by kernel smoothness, if the kernel is at least Lipschitz, and produces quasi-uniform point sets~\cite{de2005near,jeong2025convergence,santin6convergence}. 

Returning to general kernels, the above perspective is to find node sets  that works for all $G(t)$ in the range of the kernel $K$. An alternative perspective, which is closely related to practical applications, is given a specific $G = K\rho_0$ for some unknown $\rho_0$, find $\rho$ such that $G = K\rho$ and $\rho$ is as sparse as possible. In particular, such $\rho$ should have the property that $\|\rho\|_0$ should be (much) smaller than the size of the node sets $\{\omega_1,\cdots,\omega_k\}$ chosen by the cross approximation.

One route to achieve this is to take advantage of the structure of $K$. Without going into too much detail, this problem could be transformed into a rational approximation or sum-of-exponential approximation problem, so alternative methods (such as AAA algorithm) have been very useful when $K$ is an integration kernel in one dimension. Such routes are not practical for now when $K$ is an integration kernel in higher dimensions (say 2d, 3d).

\medskip
\begin{update}
Reference
\begin{quote}
Venkata Siddharth Pendyala, Local-to-Global Convergence of Greedy Cross Approximation for Totally Positive Kernels, 2026. \url{https://doi.org/10.5281/zenodo.21863274}
\end{quote}
claims to have solved Problem~\ref{problem:gecp}.
GECP applied to the (continuous) fermionic kernel above satisfies
\[
k=O\left(\log(\Lambda)\log(1/\varepsilon)\right)
\quad\Longrightarrow\quad
\|K-\widehat K\|_\infty\leq\varepsilon.
\]
Thus GECP attains, up to universal constants, the rank scaling previously
established for the existence of a low-rank approximation. The same bound also
holds for GECP applied to any finite matrix obtained by sampling $K$ at
distinct time and frequency nodes. Related progress toward Problem 4.2 was also obtained in
\begin{quote}
Marc Aurèle Gilles, Convergence rates for pivoted QR and LU, 
arXiv:2607.26863, 2026.
\end{quote}
\end{update}

\subsubsection{Orthogonal rows} 
\source{Anil Damle and Daniel Kressner}

Another common structural assumption is to assume that we want to compute a RRQR of $Q^T$ where $Q \in \mathbb R^{n\times k}$ has orthonormal columns ($k< n$). Of particular interest is the case where we want to pick $k$ columns of $Q^T.$  Addressing this task is important in, e.g., model reduction~\cite{drmac2016new} and certain algorithms for low-rank matrix approximation~\cite{golub1976rank}. 

It is known that there exists $I$ such that $\|Q(I,:)^{-1}\|_2 \le \sqrt{k(n-k+1)}$~\cite{hong1992rank} and there is a fairly efficient $O(nk^2)$ algorithm that satisfies this bound~\cite{Osinsky2025close}.    

\begin{problem} \label{problem43}
Does the QRCP algorithm typically used in practice satisfy the above, or a similar, bound.
\end{problem}
Existing upper bounds grow exponentialy in $k$~\cite{gu1996efficient}, though the known worst case examples do not have orthonormal rows. 

\medskip
\begin{update}
    Theorem 2 in \begin{quote}
        Leheng Chen, Zihao Liu, Wanyi He and Bin Dong, Iteris: Agentic Research Loops for Computational Mathematics, arXiv:2606.02484, 2026
    \end{quote}
    claims to give a negative answer to Problem~\ref{problem43}. It constructs a family of matrices $Q$ showing that, in fact, no polynomial bound of the form
    \[
    \|Q(I,:)^{-1}\|_2 \le C k^\alpha \sqrt{k(n-k+1)},
    \]
    with some constant $C$ and exponent $\alpha > 0$, can hold.
    
    On the positive side, reference
    \begin{quote}
    Anil Damle, Computing Strong Rank-Revealing Factorizations for Matrices with Orthonormal Rows, arXiv:2607.13532, 2026
    \end{quote}
shows that replacing Golub--Businger pivoting by Bischof--Stewart pivoting yields
\[
\|Q(I,:)^{-1}\|_2 \leq \sqrt{k(n-k)+1}.
\]
Thus, although the answer for standard QRCP is negative,
a closely related greedy pivoting strategy achieves the desired bound.
\end{update}

\subsubsection{Incoherence and random row subsampling}
\source{Paul Beckman}

Let $A \in \C^{n \times m}$ with $n \gg m$ be a matrix whose entries are given by a (generally non-symmetric) kernel function $A_{ij} = f(x_i, \omega_j)$. Suppose $A$ can be accurately approximated using $k$ columns with indices $J$, so that $A \approx A(:,J)T$.  

\begin{problem}
Consider selecting $r = \mathcal{O}(k)$ uniformly random row indices $I$ and then choosing column indices $\tilde{J}$ so that $A(I,:) \approx A(I,\tilde{J}) \tilde{T}$.
What are sufficient conditions on $r$ and on the kernel $f$ which guarantee that
\begin{equation}
    \|\Pi_{\tilde{J}} A - A\|_F \leq c(n, m, k) \|\Pi_J A - A\|_F 
\end{equation}
with high probability, where $c(n, m, k)$ grows at most logarithmically in its arguments? In other words, when does column selection on a random subset of the rows perform similarly to column selection on the full matrix?
\end{problem}

\motivation{When constructing compressed matrix representations of oscillatory operators such as Fourier integral operators $f(x, \omega) = e^{2\pi i \Phi(x,\omega)}$, one must compress interactions between narrow spectral bands and large physical domains. This exactly corresponds to computing low-rank approximations to matrices $A \in \C^{n \times m}$ with $n \gg m$ as described above. Computing a rank-$k$ interpolative decomposition of $A$ directly would require $\mathcal{O}(nmk)$ operations. Conversely, the algorithm described above first randomly subsamples $\abs{I}$ rows, then computes a ``good'' set of column indices $\tilde{J}$ of the $\abs{I} \times m$ matrix $A(I,:)$, which requires only $\mathcal{O}(\abs{I}mk)$ operations. This approach was introduced and showed favorable performance in the high-frequency scattering literature~\cite{engquist2009fast, pang2020interpolative}, in which one generally takes e.g. $\abs{I} = 5k$, yielding an $\mathcal{O}(k^2m)$ \textit{sublinear} column selection method. Existing error bounds on the subsampled column selection algorithm assume there exists a factorization of $A$ whose factors have suitable sparsity or incoherence properties~\cite{chiu2013sublinear,cortinovis2025sublinear}. However, such an assumption is computationally prohibitive to verify. Therefore we desire a direct, efficiently verifiable sufficient condition on the kernel $f$ under which the randomized algorithm performs well.}

\subsubsection{Block-column subset selection for rank-structured matrices}

\source{Diana Halikias}

Rank-structured matrices (e.g., $\mathcal H, \mathcal H^2$, HODLR, HSS, BLR, etc.) are useful structured approximations that come with a set of fast and convenient linear algebra routines. The ``partition'' of a rank-structured matrix determines the indices of subblocks of the matrix that are low-rank. We can think of the set of practical partitions in a problem as a set $\mathcal P$ that is finite, but combinatorially large (depending on parameters such as the rank $k$, number of low-rank blocks $B$, the maximum block size $M$, and so on). For notational simplicity, we fold these parameters into a given partition $P \in \mathcal P$. Given a matrix $A\in \R^{n \times n}$, we denote by $A_P$ the best rank-structured approximation to $A$ among all matrices with partition $P$. This operation is given by truncating to the best rank-$k$ approximation of each block that is designated as low-rank in $P$.

One can naturally extend the problem of column subset selection to block-column subset selection in each block-row of a matrix partitioned to its finest level. In particular, we can ask how to obtain a quasi-optimal approximation to $A$ among all hierarchical matrices with partitions in $\mathcal P$. Let $P^\ast = \text{argmin}_{P \in \mathcal P} \|A _{P} - A \|_F $ be the optimal partition choice.   Then, the goal is to find a partition $P \in \mathcal P$ satisfying 

\begin{equation}\label{eq:block_column_subset_selection}
    \| A_P- A \|_F \leq c(P^\ast)  \|A_{P^\ast} - A \|_F.
\end{equation}

This question is closely related to column subset selection, as there are too many blocks to feasibly test for low-rank structure. One can hope to use an adaptive/greedy method for column subset selection applied to several columns at a time. 

In partial differential equation (PDE) learning, a matrix $A$ can arise as a discretization of a Green's function (associated with a PDE operator), i.e., $A_{ij} \approx G(x_i,x_j)$. Hierarchical matrix formats exploit the fact that, for certain classes of operators, the Green's function is approximately separable on product domains $X \times Y$ away from its singular support. For example, for
elliptic operators, $G(x,y)$ is smooth for $x \neq y$ and admits accurate local expansions on well-separated sets, implying that off-diagonal blocks of $A$ have rapidly decaying singular values when the index sets are geometrically separated. In contrast, for hyperbolic or strongly advection-dominated problems, the Green's function is concentrated along characteristic manifolds, so that blocks that are well separated in Euclidean distance may nevertheless exhibit high numerical rank. In such settings, a successful partition must adapt to the anisotropic transport geometry of the operator rather than rely on a fixed, distance-based hierarchy. Consequently, the choice of partition $P$ should be viewed as encoding structural assumptions about the microlocal or transport properties of the Green's function that ensure that certain matrix blocks compressible.

\begin{problem}
For what  classes of linear PDEs (e.g., hyperbolic, elliptic, parabolic), does there exist a practical, efficient algorithm that produces a quasi-optimal approximation to the set of solution operators as in~\ref{eq:block_column_subset_selection}? What is the optimal dependence of $c$ on the partition $P$'s parameters?
\end{problem}
\motivation{Rank-structured approximations achieve  success in practice, particularly for compressing the solution operators of elliptic PDEs. If we wish to extend this methods to harder problems and other settings, it is necessary to use an adaptive or greedy method for choosing the partition, rather than impose a one-size-fits-all structure. For example, the Green's function of a hyperbolic PDE has characteristic curves that disrupt the low-rank structure of a conventional partition choice~\cite{wang2023operatorlearninghyperbolicpartial}. There is already evidence that a tailored representation can be more accurate for time-dependent PDEs~\cite{masseirobolkressner2022} and tensor formats~\cite{ehrlacher2021}, however we lack theoretical guarantees for  general settings where these algorithms can work well. Moreover,  for matrices arising in data science, where row and column indexing lacks a physical meaning, these results may also enable us to find a suitable permutation of the rows and columns to obtain a matrix that is closer to rank-structured.}

\subsubsection{Inverses of Laplacians and related matrices}
\source{Mark Fornace}

{ \newcommand{\I}{\mathcal{I}}

A new (and likely important) motivation for column subset selection comes from a problem of Markov chain compression.
In this problem, to summarize and simplify, one is presented with a large graph Laplacian $L$, and one wants to compute a reduced model in an interpolative fashion, yielding a small graph Laplacian which approximates the original one well.
In order to guarantee accuracy with respect to the long-timescale dynamics of the respective random walk, it suffices (\cite{Fornace2025-approximation}) to bound the accuracy of the Nystr\"om approximation in spectral or nuclear norm:
\begin{equation}
    \lVert K - K_{:,\I} K_{\I,\I}^{-1} K_{\I,:} \rVert_{\{2,*\}} \label{eq:nystrom-error}
\end{equation}
with respect to chosen subset $\I$ while taking $K = (L + \gamma \mathbf{I})^{-1}$ with $\gamma > 0$.
(Strictly, \cite{Fornace2025-approximation} takes this objective in the limit $\gamma \rightarrow 0^+$.)
We will focus on the nuclear error for simplicity and ease of derivation.
Note that the nuclear Nystr\"om error in (\ref{eq:nystrom-error}) may be exactly equated to the CX Frobenius error in (\ref{eq:CSSP}) by taking $A = K^{1/2}$, for instance (e.g., \cite{Fornace2024-column-01}).

In \cite{Fornace2024-column-01} it is proved that the nuclear Nystr\"om error for inverse Laplacians, taking care of null space issues, is a submodular function in $\I$ if one excludes the empty set from consideration.
Roughly speaking, this implies that the worst-case nuclear approximation error obtained by choosing $k$ columns greedily is bounded to the one obtained by choosing $s$ columns optimally by a factor decaying like $e^{-k/s}$ where $k\geq s$.
This question may be extended beyond Laplacians to consider positive-definite matirces that are either symmetric diagonally dominant (SDD) matrices or SDD M-matrices (SDDM).

\begin{problem} \label{problem46}
    Prove or disprove the submodularity of the nuclear Nystr\"om error (\ref{eq:nystrom-error}) when $L$ is assumed, in contrast to the above, to be (a) SDDM and positive-definite or (b) SDD and positive-definite.
\end{problem}
It is expected that the extension (a) is simple to prove, while (b) is less clear.
For the latter, on the other hand, it may be possible to use existing tricks mapping SDD systems to Laplacians to establish satisfactory algorithms and bounds.
It would also be interesting to consider this same question for less closely related classes of matrices.}

\medskip
\begin{update}
Reference
\begin{quote}
Matthew J. Colbrook, Nyström Error Beyond M-Matrices: A Minimal Diagonally Dominant Obstruction, arXiv:2607.19282, 2026
\end{quote}
claims to solve Problem~\ref{problem46}. It  provides a positive answer for (a). It provides a 
 negative answer for (b) (positive definite SDD), by constructing a $3\times 3$ parametrized counterexample.
\end{update}

\subsection{Clarifying relations between volume sampling and optimal subset selection}

\scribe{Mark Fornace}

{
\newcommand{\I}{\mathcal{I}}
\newcommand{\Tr}[1]{\mathrm{Tr}\left[#1\right]}

\begin{problem}
    Consider $n$ positive real numbers $\lambda$, and consider the set of all $n \times n$ orthonormal matrices $S$.
    Prove a tightness bound(s) between the (1) trace error of the optimal column subset selection of the worst case matrix $K$ possessing eigenvalues $\lambda$:
    \begin{equation}
    x_k (\lambda) := \max_{V \in S} \min_{\I \subset [n], |\I| = k} \Tr{K - K_{:,\I} K_{\I, \I}^{-1}  K_{\I,:} \biggr\rvert_{K=V^\top \mathrm{diag}(\lambda) V}}
    \end{equation}
    and (2) the objective yielded by volume sampling any such $K$:
    \begin{equation}
      y_k (\lambda) := \left(\mathbb{E}_{\I \sim k-\mathrm{DPP}(K)} \Tr{K- K_{:,\I} K_{\I, \I}^{-1}  K_{\I, :}}\right) \bigr\rvert_{K=\mathrm{diag}(\lambda)} = (k+1) \frac{e_{k+1} (\lambda)}{e_k(\lambda)}
    \end{equation}
    which is dependent on $\lambda$ alone in terms of the elementary symmetric polynomials $e_k$\cite{Guruswami2012-optimal}.
\end{problem}

\motivation{
A tightness bound, especially a simple one, could provide useful intuition as to the limits of any column subset selection bounds based on this objective and the input matrix spectrum alone.
It would also provide stronger interpretation of the bound between (1) greedy optimization of this objective and volume sampling (2) established in \cite{Fornace2024-column-01}.

For an easy edge case, consider symmetric positive definite $K$, $L = K^{-1}$, and $k = n-1$. Denote the single element $j$ (the unselected column) such that $\{j\} = [n] \setminus \I$.
Then:
\begin{equation}
  \Tr{K - K_{:,\I} K_{\I,\I}^{-1} K_{\I,:}} = L_{j,j}^{-1}
\end{equation}
by Schur complements.
For any positive real numbers $\lambda$, there exists a symmetric matrix $K$ with eigenvalues $\lambda$ such that the elements of the diagonal of $L$ are all equal, by the Schur-Horn theorem.
Therefore $x_{n-1} (\lambda) = y_{n-1} (\lambda)$ for any such $\lambda$.
However, what happens for other $k$ is unclear, although limited empirical investigation indicates that the two quantities are usually quite close.
Closely related is the question of what the variance (or other higher-order statistics) of volume sampling is.
}

}

\section{Sketching}
\label{sec:sketch}

In modern scientific computing and data science, we frequently encounter matrices of tremendous scale.
To handle such problems, \emph{sketching} techniques have proven immensely powerful.
These methods apply a linear map, typically randomly generated, to reduce the size of a large linear algebra problem, after which standard direct techniques can be used to process the data.

Sketching has proven an immensely useful paradigm in computational linear algebra.
Yet basic questions remain about the theory of sketching, and questions remain about which sketching maps to use in practice.

\subsection{Injections vs.\ embeddings}
\scribe{Raphael A.\ Meyer}

What properties must a sketching matrix have to be effective?
There are several ways of quantifying the effectiveness of a (random) matrix $\Omega$ as a dimensionality reduction map:

\begin{definition}[Subspace embeddings and injections]
    Let $\Omega \in \mathbb{R}^{n \times k}$ be a random matrix, and fix an $r$-dimensional subspace $\mathcal{V} \subseteq \mathbb{R}^n$.
    \begin{itemize}
        \item The matrix $\Omega$ is called a \emph{subspace embedding} for $\mathcal{V}$ with \emph{injectivity} $\alpha$ and \emph{dilation} $\beta$ if 
        \begin{equation*}
            \alpha \| x \|_2^2 \le \|\Omega^\top x\|_2^2 \le \beta \| x \|_2^2 \quad \text{for every } x \in \mathcal{V}.
        \end{equation*}
        \item The matrix $\Omega$ is called a \emph{subspace injection} for $\mathcal{V}$ if it satisfies two conditions:
        \begin{enumerate}
            \item \textbf{\textit{Isotropy:}} For every vector $x \in \mathbb{R}^n$, it holds that $\mathbb{E} \|\Omega^\top x\|_2^2 = \|x\|_2^2$.
            \item \textbf{\textit{Injectivity:}} It holds that 
            \begin{equation*}
            \alpha \| x \|_2^2 \le \|\Omega^\top x\|_2^2 \quad \text{for every } x \in \mathcal{V}.
        \end{equation*}
        \end{enumerate}
    \end{itemize}
    A matrix $\Omega$ is called an \emph{oblivious} subspace embedding (resp.\ injection) with failure probability $\delta$ if, for any fixed $r$-dimensional subspace $\mathcal{V} \subseteq \mathbb{R}^n$, the matrix $\Omega$ is an subspace embedding (resp.\ injection) for $\mathcal{V}$ with at least $1-\delta$ probability.
    If $\delta$ is not specified, we take it to be some arbitrary small value, say $\delta = 0.01$.
\end{definition}

The subspace embedding property, essentially due to Sarl\'os \cite{Sar06}, and many linear algebra algorithms have been analyzed under the assumption that $\Omega$ is a subspace embedding; see surveys \cite{Woo14b,KT24}.
However, there are reasons to believe that the subspace embedding property is too restrictive to capture the behavior of certain sketching matrices $\Omega$.
The subspace injection property was investigated by Cama\~no, Epperly, Meyer, and Tropp \cite{CEMT25a} as a weaker condition that still suffices to justify that many randomized linear algebra algorithms produce solutions that are accurate to constant-factor accuracy.
Conditions similar to the subspace injection condition were also investigated in earlier works (e.g., \cite{DMMS11}).
The terms oblivious subspace embedding and oblivious subspace injection are abbreviated OSE and OSI.

For many sketching matrices in common, we have OSE guarantees.
However, for certain highly structured embeddings, only OSI guarantees are available (in some parameter regimes). 
Therefore, it is natural to inquire as the minimal requirements on $\Omega$ to justify any particular algorithm.
\begin{quote}
    \textbf{Meta-problem:} For a given algorithm that uses a sketching matrix $\Omega$, what is the weakest assumption we can make on $\Omega$ that ensures a desired level of accuracy?
\end{quote}
Examples of conditions on the sketching matrix are a symmetric OSE guarantee ($\alpha = 1-\varepsilon$, $\beta = 1+\varepsilon$), an asymmetric OSE guarantee ($\alpha$, $\beta$ arbitrary), or an OSI guarantee (possibly with side guarantees, such as on approximate matrix multiplication).
Examples of algorithms include the sketch-and-solve and sketch-and-precondition algorithms for linear least-squares, and the randomized SVD for low-rank approximation.
See \cite{KT24} for introductions to these algorithms.

Cama\~no et al.\ \cite{CEMT25a} showed that OSIs with constant injectivity sufficed to solve least squares regression and low-rank approximation to constant factor accuracy.
A first suite of open problems asks if a higher quality OSI suffices for a relative-error guarantee:

\begin{problem}[Sketch-and-solve]
    \label{prob:sketch-and-solve}
    Let \(A \in \mathbb R^{n \times d}\) be a full-rank matrix, \(B\in\mathbb R^{n \times p}\) be a matrix, and let \(\Omega \in \mathbb R^{n \times k}\) be an oblivious subspace injection with dimension $d$ and injectivity \(1-\varepsilon\).
    Let \(\tilde X = (\Omega^\top A)^+(\Omega^\top B) \in \mathbb R^{d \times p}\) be the sketch-and-solve approximation to the ordinary least-squares solution $X = A^+B$.
    (Here, ${}^+$ denotes the Moore--Penrose pseudoinverse.)
    Prove or disprove: Does it necessarily hold that
    \[
        \|A \tilde X - B\|_{\rm F}
        \leq (1+\mathcal{O}(\varepsilon)) \min_{X \in \mathbb R^{d \times p}}\|AX - B\|_{\rm F}?
    \]
\end{problem}
\begin{update}
Section 3.1 of
\begin{quote}
Alex Townsend and Christopher Wang, Oblivious Subspace Injection Is Not Enough for Relative Error, arXiv:2604.10215v2, 2026.
\end{quote}
claims to give a negative answer to Problem~\ref{prob:sketch-and-solve}.  
It constructs a linear least-squares problem (with a $2\times 1$ matrix $A$) and an OSI with injectivity $1-\varepsilon$ for which the sketch-and-solve least-squares estimator does not imply a $1+\mathcal O(\varepsilon)$ guarantee with failure controlled solely by the OSI failure parameter
The broader problem of identifying the weakest convenient additional condition remains open.
\end{update}
\medskip

Likewise, one can ask this question for low-rank approximation:
\begin{problem}[Randomized SVD]
    \label{prob:rsvd}
    Let \(A \in \mathbb R^{n \times d}\) be a matrix, fix $r > 0$, and let \(\Omega \in \mathbb R^{n \times k}\) be an oblivious subspace injection with dimension $r$ and injectivity \(1-\varepsilon\).
    Let $\tilde{A} = (A\Omega)(A\Omega)^+A$ denote the randomized SVD approximation to $A$, and let $A_r$ be the best rank-$r$ approximation to $A$ with respect to the Frobenius norm.
    Prove or disprove: Does it necessarily hold that
    \[
        \|A - \tilde{A}\|_{\rm F}
        \leq (1+\mathcal{O}(\varepsilon)) \|A-A_r\|_{\rm F}?
    \]
\end{problem}
\begin{update}
Section 3.2 of
\begin{quote}
Alex Townsend and Christopher Wang, Oblivious Subspace Injection Is Not Enough for Relative Error, arXiv:2604.10215v2, 2026
\end{quote}
claims to give a negative answer to Problem~\ref{prob:rsvd} by constructing a 
$2\times 2$ counterexample.
\end{update}
\medskip

Christopher Musco has an analysis showing that for the least-squares problem (Problem \ref{prob:sketch-and-solve}), a sketching matrix \(\Omega\) that is an OSI with injectivity \(1-\varepsilon\) and with failure probability \(\delta=\mathcal O(\varepsilon)\) solves Problem \ref{prob:sketch-and-solve} when \(p = 1\).
As a result, if it is shown that Problems \ref{prob:sketch-and-solve} or \ref{prob:rsvd} are not solvable using an OSI that has constant failure probability, we ask if a sufficiently small failure probability (say, \(\delta = \poly(\varepsilon)\)) does suffice.

Lastly, we can ask about generalizing the OSI condition to $\ell_p$ regression problems.

\begin{problem}[Regression in the $\ell_p$ norm]
\label{problem53}
    Consider the $\ell_p$ linear regression problem
    \begin{equation*}
        \operatorname{minimize}_{x \in\mathbb{R}^d} \|Ax - b\|_{\ell_p}
    \end{equation*}
    for a matrix \(A \in \mathbb R^{n \times d}\) and a vector \(b\in\mathbb R^{n}\).
    For a sketching matrix \(\Omega \in \mathbb R^{n \times k}\), the sketch-and-solve solution is \[\tilde x \in \operatorname{argmin}_{x \in\mathbb{R}^d} \|(\Omega^\top Ax)-(\Omega^\top b)\|_{\ell_p}.\]
    Formulate a generalization of the OSI property to the $\ell_p$ regression problem which ensures that this procedure succeeds to constant-factor accuracy 
    \[
        \|A \tilde x - b \|_{\ell_p}
        \leq \mathrm{const}\cdot \min_{x \in \mathbb R^d}\|Ax - b\|_{\ell_p}.
    \]
\end{problem}
\begin{update}
Section 5 of
\begin{quote}
Alex Townsend and Christopher Wang, Oblivious Subspace Injection Is Not Enough for Relative Error, arXiv:2604.10215v2, 2026
\end{quote}
claims to solve Problem~\ref{problem53}
by providing a $p$-isotropic $\ell_p$-OSI property and 
and proving that it yields a constant-factor sketch-and-solve guarantee for $1 \le p < \infty$.
\end{update}

\subsection{Sparse dimensionality reduction maps} 
\scribe{Ethan N.\ Epperly, Eliza Rebrova, and Micha\l{} Derezi\'nski}

Sparse dimensionality reduction maps are increasingly establishing themselves as a fundamental tool in matrix computations (see, e.g., \cite{TYUC19,derezinski2021sparse,DM23,MBM+23,chenakkod2024optimal,CNR+25,CEMT25a,EMN25,FG25}).
In some computing environments, they often achieve significant speedups over other types of dimensionality reduction maps \cite{TYUC19,DM23,CEMT25a,CNR+25,Epp25a}.

There are several constructions of sparse dimensionality reduction maps, some of which behave similarly to one another and others which have vastly divergent behaviors.
To focus discussion, we consider a single model, originally introduced by Kane \& Nelson \cite{KN12}:

\begin{definition}[SparseStack \protect{\cite{KN12}}]
    Fix parameters $\zeta,b \ge 1$ and set $k \coloneqq b\zeta$.
    A SparseStack random matrix 
    \begin{equation*}
        \Omega = \frac{1}{\sqrt{\zeta}}\begin{bmatrix}
            \varrho_{11} e_{s_{11}}^* & \cdots & \varrho_{1\zeta} e_{s_{1\zeta}}^* \\
            \varrho_{21} e_{s_{21}}^* & \cdots & \varrho_{2\zeta} e_{s_{1\zeta}}^* \\
            \vdots & \ddots & \vdots \\
            \varrho_{n1} e_{s_{n1}}^* & \cdots & \varrho_{n\zeta} e_{s_{n\zeta}}^*
        \end{bmatrix} \in \mathbb{R}^{n \times k}.
    \end{equation*}
    Here, $\varrho_{ij}$ are iid Rademacher random variables, and $s_{ij}$ are iid uniform random indices on $\{1,\ldots,b\}$.
\end{definition}

For over a decade, researchers in different communities have debated: How should we pick the parameters $\zeta$ and $k$ to ensure that SparseStacks (and related constructions) are effective as sketching matrices?
Traditionally, this question has been studied using the oblivious subspace embedding condition; more recently, researchers have considered the weaker oblivious subspace \emph{injection} problem.
This question is more than a mathematical curiosity: Practitioners need guidance about how to choose parameters to ensure reliable performance of sketching algorithms.

The optimal embedding properties of sparse random embeddings are the subject of a conjecture of Nelson \& Nguyen \cite{NN13a}. 

\begin{problem}[Nelson--Nguyen conjecture: Subspace embedding properties of sparse random matrices] \label{problem54}
    Prove or disprove the Nelson--Nguyen conjecture: A SparseStack embedding with parameters $\zeta = \mathcal{O}((\log r) / \varepsilon)$ and $k = \mathcal{O}(r / \varepsilon^2)$ is an oblivious subspace embedding with parameters $\alpha = 1 - \varepsilon$ and $\beta = 1+\varepsilon$.
\end{problem}

Their conjecture remains open, though we are painfully close to resolving it.
Existing guarantees offer different trade-offs between $\zeta$ and $k$ that approach but do not match the conjecture. Most notably, \cite{cohen2016nearly} matched the conjectured sparsity $\xi=O((\log r)/\varepsilon)$ but with sub-optimal $k=O(r(\log r)/\varepsilon^2)$, while \cite{chenakkod2025optimal} matched the conjectured dimension $k=O(r/\varepsilon^2)$, but with sparsity $\xi = O((\log r)^3 + (\log r)^2/\varepsilon)$. Finally, \cite{chenakkod2026optimal} matched both parameters up to sub-polylogarithmic factors, i.e., $\xi = O((\log r)^{1+o(1)}/\varepsilon)$ and $k=O((r/\varepsilon^2)(\log r)^{o(1)})$, where the exponent $o(1)$ scales as $1/(\log\log\log\log r)$.

\medskip
\begin{update}
The reference
\begin{quote}
    Joel A. Tropp, Comparison Theorems for the Extreme Eigenvalues
    of a Random Symmetric Matrix, arXiv:2603.04365v4, 2026
\end{quote}
claims significant progress on Problem~\ref{problem54}. For a complex
SparseStack variant, it proves the OSI property
with injectivity $1-\varepsilon$ and, for constant failure probability,
$\zeta = \mathcal{O}((\log r) / \varepsilon)$ and $k = \mathcal{O}(r / \varepsilon^2)$.
Thus, the lower-distortion bound is obtained at the parameter scaling
conjectured in Problem~\ref{problem54}. The corresponding
upper-distortion bound remains open.

The reference
\begin{quote}
Diar Heidary, A Sparse-Fock Proof of the Upper Edge for Fully Independent SparseStack,
Version 1.5, 2026. \url{https://github.com/DiarHaidary/nelson-nguyen-sparse-fock}
\end{quote}
claims to settle Problem~\ref{problem54} in the affirmative, proving the full two-sided OSE property with $\zeta = O((\log r)/\varepsilon)$ and $k = O(r/\varepsilon^2)$ for constant failure probability, matching the parameter scaling conjectured by Nelson and Nguyen.
\end{update}
\medskip

The full version of the Nelson--Nguyen conjecture concerns the dependence of $k$ and $\zeta$ on both $r$ and $\varepsilon$, but the following weaker question is also interesting: is a SparseStack with parameters $k = \mathcal{O}(r)$ and $\zeta = \mathcal{O}(\log r)$  an oblivious subspace embedding where $\alpha,\beta$ are any two universal constants?

Practitioners have pushed farther than Nelson \& Nguyen's recommended parameter settings $k = \mathcal{O}(r)$ and $\zeta = \mathcal{O}(\log r)$, investigating the use of sparse dimensionality reduction maps with sparsity parameters $\zeta$ as low as $2$, with strong empirical evidence of reliability with $\zeta \ge 4$.
However, Nelson \& Nguyen proved that no sparse embedding can be an oblivious subspace \emph{embedding} with constant parameters $\alpha$ and $\beta$ cannot have sparsity $\zeta = o(\log r / \log \log r)$ \cite{NN14}.
For this reason, the subspace \emph{injection}  condition provides a path forward to understanding the empirical success of sparse dimensionality maps with constant $\zeta = \mathcal{O}(1)$.

\begin{problem}[Subspace injection properties of sparse random matrices] \label{problem55}
    Assume the condition  $r\le \mathcal{O}(n/\log n)$ relating the subspace dimension $r$ and the ambient dimension $n$.
    Prove or disprove: A SparseStack embedding with parameters $\zeta = \mathcal{O}(1)$ and $k = \mathcal{O}(r)$ is an oblivious subspace injection with parameter $\alpha = \mathrm{const}$.
    More specifically, do parameters $\zeta = 4$ and $k = 2r$ suffice to achieve the OSI property (for all sufficiently large $r$)?
\end{problem}

Tropp \cite{Tro25} nearly achieved this result, by showing that $\zeta = \mathcal{O}((\log r)/\varepsilon^2)$ and $k = \mathcal{O}(r/\varepsilon^2)$ achieve to produce an OSI with $\alpha = 1-\varepsilon$.
Konstantin Tikhomirov presented an argument to a small group at the Simons workshop that suggests the condition $r\le \mathcal{O}(n/\log n)$ may be necessary.

\medskip
\begin{update}
Reference
\begin{quote}
Han Huang, Mark Rudelson and Konstantin Tikhomirov,
Well-invertible column subsets of sparse matrices are rare, arXiv:2607.05384v2, 2026
\end{quote}
claims to answer Problem~\ref{problem55} negatively in a broader setting. Corollary 1.6 states that for every fixed $\zeta \ge 2$, every fixed $\varepsilon, \alpha >0$, sufficiently large block size, and $n\ge Ck$, a fixed coordinate subspace of dimension at least $\varepsilon k$ has minimum sketched singular value below $\alpha$ with probability at least $0.99$. In particular, this rules out $\zeta = 4$, $k = 2r$, and constant $\alpha$.
\end{update}

\subsection{Subsampled randomized trigonometric transforms}
\scribe{Ethan N.\ Epperly}

Another popular class of dimensionality reduction maps are based on fast transforms. 
We make the following definition:

\begin{definition}[Subsampled randomized Hadamard transform]
    A (traditional) subsampled randomized Hadamard transform (SRHT) $\Omega \in \mathbb{R}^{n\times k}$ is a random matrix of the form 
    \begin{equation*}
        \Omega = \sqrt{\frac{n}{k}} DFS,
    \end{equation*}
    where $D$ is a diagonal matrix with iid Rademacher entries, $F$ is the Walsh--Hadamard matrix, and $S$ subsamples down to $k$ uniformly random positions.
    We may also define a \emph{rerandomzied} SRHT as 
    \begin{equation*}
        \Omega = \sqrt{\frac{n}{k}} D_1FD_2FS.
    \end{equation*}
    Here, $D_1$ and $D_2$ are independent diagonal matrices with iid Rademacher entries.
\end{definition}

Traditional SRHTs, and variants using other types of fast transform, are widely used in randomized linear algebra. 
The following theorem, due to Tropp \cite{Tro11}, analyzes the behavior of traditional SRHTs.

\begin{theorem}[Subsampled randomized Hadamard transform]
    A traditional SRHT is a subspace embedding with $\alpha,\beta = \mathrm{const}$ with embedding dimension $\mathcal{O}((r + \log n) \log r)$.
    Moreover, there is an $r$-dimensional subspace $\mathcal{V} \subseteq \mathbb{R}^{r^2}$ such that a traditional SRHT is \emph{not} a subspace injection if $k = o(r \log r)$.
\end{theorem}

Assume $r\ge \log n$, this result shows that $r \log r$ is the correct scaling for the embedding dimension $k$ as a function of the subspace dimension $r$ for the traditional SRHT.
However, researchers in practice often want to use an embedding with output dimension $k = O(r)$.
Empirically, rerandomized SRHTs appear to observe this scaling.

\begin{problem}[Rerandomized subsampled trigonometric transforms]
    Assume that $n = \Omega(\log r)$.
    Prove or disprove: The rerandomized SRHT of dimension $k = \mathcal{O}(r/\varepsilon^2)$ is an oblivious subspace embedding with $\alpha = 1-\varepsilon$ and $\beta = 1+\varepsilon$.
\end{problem}

As with the previous section, weaker versions of this questions are also interesting, such as where we allow the parameters $\alpha$ and $\beta$ to be arbitrary constants or where we only require the embedding to be an OSI.

During the workshop, Sarl\'os communicated to us that a similar construction to the rerandomized SRHT is used in the field of locality sensitive hashing (LSH) \cite{AIL+15}.
In this context, the authors report that two rounds of randomization, as in the rerandomized SRHT, are empirically not sufficient for LSH purposes, but three rounds do suffice.
Unfortunately, they do not provide an explicit example on which two levels of randomization fails for LSH.
Some guarantees for random orthogonal matrices of the form $D_1FD_2FD_3F$ are provided in \cite{CFG+16}.

\section{Tensor, quantum, and other problems}
\label{sec:other}
\subsection{Optimal bounds for Tensor Train decomposition}
\scribe{Chris Camaño, Edgar Solomonik, and Jess Williams}
\source{This question was initially proposed by Mehrdad Ghadiri.}

Given a tensor $\mathcal{X}\in\mathbb{R}^{d_1\times\cdots \times d_n}$, a tree tensor decomposition of $\mathcal{X}$ is a contraction of a tensor network of $m$ tensors, with $n-1\leq m \leq 2n-1$, such that each tensor is assigned at most 1 index (mode) of $\mathcal{X}$ and each leaf in the tensor network is assigned exactly 1 index of $\mathcal{X}$.
Moreover, if non-leaf nodes exist in the tree, there must be a choice of root node such that each internal node has at least two descendants.
For example, a tensor train decomposition is a flat tree tensor network where each tensor factor is assigned one index of $\mathcal{X}$ (so $m=n$), and a Tucker decomposition is a star graph (height 1 tree) with $m=n+1$.
In further detail, a tensor train decomposition with ranks $\chi_1,\ldots,\chi_{n-1}$, corresponds to the factorization,
$$
\mathcal{X}(i_1,\ldots,i_n)
=\sum_{\alpha_1=1}^{\chi_1}\cdots\sum_{\alpha_{n-1}=1}^{\chi_{n-1}}
  \mathcal{G}_1(i_1,\alpha_1)\,
  \Bigl(\prod_{k=2}^{n-1} \mathcal{G}_k(\alpha_{k-1},i_k,\alpha_k)\Bigr)\,
  \mathcal{G}_n(\alpha_{n-1},i_n).
$$

\motivation{
Approximation of tensors with such decompositions is a widely applicable and well-studied task~\cite{kolda2009tensor}.
However, the best known general theoretical bounds on the approximation error of tree tensor decompositions are limited to bounds attained from approximations attained by a sequence of low-rank matrix approximations (SVDs).
For tensor train, the corresponding algorithm is referred to as TT-SVD~\cite{oseledets2011tensor}, while for Tucker decomposition, it corresponds to the sequentially truncated high-order SVD (HOSVD) algorithm.
These algorithms perform a sequence of $m-1$ SVDs (low rank matrix approximations), to recursively split up the tensor according to the structure of the tree tensor decomposition.

For an order $m$ tensor $\mathcal{A}$, suppose there exists a tree decomposition with tree $T=(V,E)$ yielding approximation $\mathcal{X}$ of a given Frobenius-norm approximation error, $\|\mathcal{A}-\mathcal{X}\|_F^2\leq \epsilon.$
The existence of $\mathcal{X}$, implies that for any edge $e\in E$, a low-rank matrix approximation, $X_e$, of a matricization of the tensor, $A_e$, (with rows and columns of the matrix corresponding to the combination of tensor indices induced by the two partitions of the tree split by the edge) exists with no greater error,
\[\|\mathcal{A}-\mathcal{X}\|_F^2\leq \epsilon \quad \Rightarrow \quad \forall e \in E, 
\|A_e-X_e\|_F^2\leq \epsilon.\]
Since, the factors obtained from any such SVD are orthogonal projections of $\mathcal{A}$, the low rank approximation error corresponding to matricizations of the factors for other edges in the tree is not increased (see~\cite{oseledets2011tensor} for a full derivation of this argument for tensor train, the generalized argument here is based on our working group discussion).
Hence, if for all edges $e$, the approximation error is of square norm $\epsilon$, the total approximation error obtain by the sequence of $m-1$ SVDs to produce the approximation $\mathcal{X}$ is at most $(m-1)\epsilon$ i.e.,
\[\forall e\in E, \|A_e-X_e\|_F^2\leq \epsilon \quad \Rightarrow \quad \|\mathcal{A}-\mathcal{X} \|_F^2\leq (m-1)\epsilon.\]

For Tucker decomposition, the above error bound has been shown to be tight for HOSVD of tensors of certain order and size~\cite{fahrbach2025tight}.
Further, we are not aware of polynomial-time algorithms for any tree decomposition (including Tucker and tensor train), that obtain a general improvement to this error bound.
Hence, we have the following open question.
}
\begin{problem}
Given any choice of $m$-node tree tensor network (e.g., tensor train or Tucker) and choice of ranks (dimensions of edges in the tensor network), let the set of representable tensors in this tensor network be $S\subset  \mathbb{R}^{d_1\times \cdots \times d_n}$.
Find a polynomial-time algorithm that takes as input any tensor $\mathcal{A}\in\mathbb{R}^{d_1\times \cdots \times d_n}$ and yields an approximation $\mathcal{X}\in\mathbb{S}$ such that
\[\|\mathcal{X}-\mathcal{A}\|_F^2 < (m-1)\min_{\mathcal{Y}\in S}\|\mathcal{Y}-\mathcal{A}\|_F^2,\]
or show that an $(m-1)$-approximation is optimal under complexity-theoretic assumptions.
\end{problem}
The problem asks for any algorithm with error that achieves a strict inequality, though ideally more substantial improvements (e.g., $o(m)$ instead of $m-1$) would be desired (indeed a $o(m)$ bound was the version of the question posed for Tucker prior to the workshop).
This problem can also be posed for hierarchical semi-separable (HSS) approximation of matrices.
The HSS approximation has a direct correspondence to approximation with low-rank tensor networks~\cite{ceruti2025low} and some of the best-known randomized HSS approximation algorithms achieve analogous error bounds in expectation~\cite{amsel2025quasioptimalhierarchicallysemiseparablematrix}.
An additional variant of the problem is to restrict the approximation problem to positive tensors or elementwise random tensors with a fixed positive mean.
Existing results have shown that contraction of positive tensor networks may be easier to approximate than general tensor networks~\cite{TN-positive-bias,TN-sign-problem}.

\subsection{Eigenvalue systems arising from Hamiltonians}

\scribe{Edgar Solomonik, discussed with Elias Jarlebring, Florian Schafer, Maryam Dehghan, Tamara Kolda, Chris Camaño}
\source{This question is based on discussion of a pre-workshop question of Elias Jarlebring on companion-type linearization and nonlinear eigenvalue problems.}

\begin{definition}
A nonlinear eigenvalue problem with eigenvector nonlinearity (NEPv)
~\cite{bai2018robust} is defined by a matrix-valued function, $A$, $A(z)\in\mathbb{C}^{n\times n}, \forall z\in\mathbb{C}^n$.
An eigenpair $(\lambda, z)$ for this problem is defined by 
\[A(z)z = \lambda z.\]
\end{definition}
\begin{problem}
Consider $n,d\in\mathbb{N}$.
Given $X\in\mathbb{C}^{n \times  n}$, for any $k\in\{1,\ldots,d\}$, we denote 
\[X^{\langle k \rangle} = I_n^{ \otimes (k-1)} \otimes X \otimes I_{n}^{\otimes (d-k)},\]
where $I_n$ is the $n\times n$ identity matrix.
Given a (Hamiltonian) matrix $H$,  
with tensor product structure,
\begin{align}
H = \sum_{i=1}^d \bigg(F^{\langle i\rangle }_i + \sum_{j=1}^d G_{i,j}^{\langle i \rangle}K_{i,j}^{\langle j \rangle}\bigg).
\label{eq:prob_bos0}
\end{align}
Assume, for all $i,j$, the matrices $F_i, G_{i,j}, K_{i,j}\in\mathbb{C}^{n\times n}$ are Hermitian and norm at most $1$.
Consider the problem,
\begin{align}
\min_{x_1,\ldots,x_d\in\mathbb{C}^n}
f(x_1,\ldots,x_d),\quad f(x_1,\ldots,x_d) =
\frac{
(x_1\otimes \cdots \otimes x_d)^HH(x_1\otimes \cdots \otimes x_d)}
{\prod_{j=1}^d x_j^Hx_j}
\label{eq:prob_bos1}
\end{align}
We can define a NEPv problem that is equivalent to \eqref{eq:prob_bos1} and of dimension $nd$,
\begin{align}
A\bigg(\begin{bmatrix}x_1 \\ \vdots \\x_d\end{bmatrix}\bigg) 
&= \begin{bmatrix} A_1 & 
& \\ & \ddots & \\ & & 
A_d\end{bmatrix}, \forall i, A_i= \bigg(\frac{1}{x_i^Hx_i}\sum_{j=1}^d x_j^Hx_j\bigg)\bigg(F_i + \sum_{j=1}^dG_{i,j}\frac{x_j^HK_{i,j}x_j}{x_j^Hx_j}\bigg).
\end{align}
Then, we have that with 
\[z = \begin{bmatrix}x_1 \\ \vdots \\x_n\end{bmatrix}, \quad
\frac{z^HA(z)z}{z^Hz}=f(x_1,\ldots,x_d).\]
Moreover, the minimizer corresponds to the solution to the NEPv defined by $A$~\cite{bai2018robust}.
Devise a computationally efficient algorithm for this NEPv or an alternative NEPv equivalent to \eqref{eq:prob_bos1}, that achieves superlinear local convergence and has iteration complexity that depends linearly on $d$.
\end{problem}
\motivation{
The standard eigenvalue problem (minimization over general vectors instead of over tensor products) for a 2-local Hamiltonian like \eqref{eq:prob_bos0} is QMA-complete~\cite{kempe2006complexity}.
The above construction is analogous (developed during the workshop and intended to provide a simple algebraic simplification) to numerical solvers for the Gross-Pitaevskii model for the time-independent bosonic Schr\"odinger equation~\cite{henning2025gross}.
A similar derivation in the fermionic case leads to the Roothaan equations and Kohn-Sham equations for density functional theory, except the eigenvector dependent nonlinear eigenproblem considers (the matrix-valued function is parameterized by) a set of eigenvectors.
The convergence guarantees and acceleration methods of such "self consistent field" procedures has been studied in literature for Hamiltonians arising in quantum chemistry, see e.g.,~\cite{yang2009convergence} and~\cite{kempe2006complexity} for a nonlinear eigenvalue problem perspective.

A general, quadratically convergent method for NEPvs such as the above is given in~\cite{bai2018robust}, but it requires the Hessian of $f(z)=z^HA(z)z$, which is $nd\times nd$.
To the best of our knowledge, there has been limited study of the complexity and convergence properties of optimization methods for such NEPvs, and the above open problem is one of many open questions that might be asked to better understand this family of methods.

\subsection{Matrix sign function}
\scribe{Elias Jarlebring and Noah Amsel}

Let $\Pi_{2^m}^*$ be the set of univariate polynomials whose corresponding matrix function can be computed with $m$ matrix-matrix multiplications and an arbitrary number of matrix additions and scalings.  
We consider the problem of determining the best such polynomial that approximates the sign function in $I:=[-1,-\delta]\cup [\delta,1],$
\begin{equation}\label{eq:sign_approx}
\epsilon_{m}^*=\min_{p\in\Pi_{2^m}^*} \max_{x\in I} |p(x)-\mathrm{sign}(x)|.
\end{equation}

\begin{problem}\label{problem:sign_main}
What is the asymptotic error $\epsilon_m^*$ as a function of $m$ and $\delta$?
\end{problem}
An important property in standard approximation theory is the equioscillation property, which is well-established for minimizers over the set of all polynomials of a given degree \cite[Thm~10.1]{Trefethen2013}.
This property states that for the degree $n$ polynomial achieving the best $L^\infty$ approximation to any given continuous function, there exists an increasing sequence of $2n$ points at which the maximum error is attained, and this error alternates in sign.
\begin{problem}
Let $p^*$ be a minimizer to the analogue \eqref{eq:sign_approx} with approximation error measured in $\ell^2$-norm. Does $p^*(x) - \mathrm{sign}(x)$ satisfy the equioscillation property?
\end{problem}
\remarks{Elias has observed empirically that equioscillation seems to hold if the distance between $p$ and $\mathrm{sign}$ is measured in the $\ell^2$ norm instead of the $\ell^{\infty}$ norm.}

Consider the class of polynomials that are represented as a composition of $T$ polynomials of the form $p_t(x) = a_t x + b_t x^3$, which we denote $\Pi_{2^{(2T)}}^{\mathrm{comp}}$.
By definition, $\Pi_{2^{2T}}^{\mathrm{comp}} \subset \Pi_{2^{2T}}^{*}$ because each $p_t(\mathbf A)$ can be computed in two matrix-matrix mulitplications.
Let
\begin{equation}\label{eq:comp_sign}
\epsilon_{2T}^{\mathrm{comp}}=\min_{p\in\Pi_{2^{2T}}^{\mathrm{comp}}} \max_{x\in I} |p(x)-\mathrm{sign}(x)|.
\end{equation}
An open problem is how much extra work must be done if we restrict ourselves to using these compositions:
\begin{problem}\label{problem:epsilon_comp}
How does $\epsilon_{2T}^{\mathrm{comp}}$ of \eqref{eq:comp_sign} compare to $\epsilon_{2T}^*$ of \eqref{eq:sign_approx}?
Alternatively, what is the smallest $T$ as a function of $m$ for which $\epsilon_{2T}^{\mathrm{comp}} \leq \epsilon_m^*$? 
\end{problem}

\motivation{
Let $\mathbf A$ be a square symmetric matrix with eigendecomposition $\mathbf V \mathbf \Lambda \mathbf V^\top$.
Define $f(\mathbf A) = \mathbf V f(\mathbf \Lambda) \mathbf V^\top$, where $f(\mathbf \Lambda)$ is the diagonal matrix with entries $f(\lambda_i)$.
We wish to approximate $\mathrm{sign}(\mathbf A)$ without taking the eigenvalue decomposition.
In fact, to be as GPU-friendly as possible, we wish to use only the operations of matrix-matrix multiplication, matrix-matrix addition, and matrix scaling, where the $O(n^3)$ cost of the former dominates.
In other words, we seek $p$ for which $p(\mathbf A) \approx \mathrm{sign}(\mathbf A)$ for all $\|\mathbf A\|_2 \leq 1$ with condition number $\leq 1/\delta$, which in turn implies that the eigenvalues are in the interval we want to approximate 
$\lambda(\mathbf A)\subset I=[-1,-\delta]\cup[\delta,1]$.

The corresponding problem for the general set of polynomials, is a classic problem in approximation theory; see  \cite{EremenkoYuditskii2007}. In this setting, the minimal degree polynomial which approximates $\mathrm{sign}$ to error $\epsilon$ in the $\ell^\infty$ norm on $I$ has degree $d=\tilde \Omega(\kappa\log(1/\epsilon))$. Using the polynomial evaluation method of \cite{doi:10.1137/0202007}, we can compute that polynomial in $m=O(\sqrt{\kappa}\log(1/\epsilon))$ matrix-matrix operations.
However, given a budget of $m$ matrix-matrix multiplications, we can hope to achieve better error than this approach due to the fact that \emph{some} degree $d$ polynomials can be computed using far fewer than $\sqrt{d}$ multiplications.
A solution to \ref{problem:sign_main} would require an improvement to this estimate.

The set $\Pi_{2^m}^*$ has been studied in \cite{jarlebring2025polynomialsetassociatedfixed}. In particular it has been shown that it is a semi-algebraic set and has a dimension $\dim(\Pi_{2^m}^*)=m^2$. Moreover,  the largest polynomial subset that is included in $\Pi_{2^m}^*$ as defined by $d_*(m)=\max\{d:\Pi_d\subset \bar{\Pi}_{2^m}^*\}$ has been characterized for $m\le 7$. As a consequence of the dimension property, one can expect that $m^2$ conditions can be satisfied, e.g., one minimizer $p^*$ might interpolate the function in $m^2$ interpolation points.

A similar problem appears in the context of electronic structure calculation. More precisely, the method class called \emph{linear scaling algorithms} involve the computation of the spectral projector of a symmetric matrix $A$, which can be computed by applying the shifted heaviside function.   For an appropriate shift $\mu$, we have $h(\mu I-A)=VV^T$, where $V\in\mathbb{R}^{n\times p}$ consists of the $p$ eigenvectors corresponding to the eigenvalues to the left of $\mu$. The linear scaling methods in this context compute approximations of the matrix heaviside function using matrix-matrix multiplications, in order to enable parallelization in a high-performance computing setting. Methods are described in \cite{Rubensson2011,Niklasson2002,DanielsScuseria1999,chen2014stable} with ideas that can be traced back to \cite{McWeeny1960}. An example of a method from this body of literature is the alternating application of the two polynomials $p_a(x)=x^2$ and $p_b(x)=x(2-x)$ - an iteration which converges to the heaviside function with $\mu=(\sqrt{5}-1)/2$. Such a construction leads to a composite polynomial, which can also be transformed to approximate the 
sign function. The full understanding 
of this method class would be insightful for 
$\epsilon_{2T}^{\rm comp}$ and Problem~\ref{problem:epsilon_comp}. 

Approximating the matrix sign function, $\mathrm{sign}(\mathbf A)$ and the related polar decomposition arises in several applications:
\begin{itemize}
\item Divide-and-conquer recursive eigensolvers,
\item Matrix approximation~\cite{higham1994matrix},
\item ``Purification'' in electronic structure calculations.
\end{itemize}

}

\bibliographystyle{alpha}
\bibliography{refs}

\newcommand{\etalchar}[1]{$^{#1}$}
\begin{thebibliography}{BEMOC21}

\bibitem[ABB{\etalchar{+}}18]{Armentano-et-al}
Diego Armentano, Carlos Beltr{\'a}n, Peter B{\"u}rgisser, Felipe Cucker, and
  Michael Shub.
\newblock A stable, polynomial-time algorithm for the eigenpair problem.
\newblock {\em J. Eur. Math. Soc. (JEMS)}, 20(6):1375--1437, 2018.

\bibitem[ACK{\etalchar{+}}25]{amsel2025quasioptimalhierarchicallysemiseparablematrix}
Noah Amsel, Tyler Chen, Feyza~Duman Keles, Diana Halikias, Cameron Musco,
  Christopher Musco, and David Persson.
\newblock Quasi-optimal hierarchically semi-separable matrix approximation,
  2025.

\bibitem[ADGJT21]{Grigori2021}
Hussam Al~Daas, Laura Grigori, Pierre Jolivet, and Pierre-Henri Tournier.
\newblock A multilevel schwarz preconditioner based on a hierarchy of robust
  coarse spaces.
\newblock {\em SIAM J.~Sci.\ Comput.}, 43(3):A1907--A1928, 2021.

\bibitem[AIL{\etalchar{+}}15]{AIL+15}
Alexandr Andoni, Piotr Indyk, Thijs Laarhoven, Ilya Razenshteyn, and Ludwig
  Schmidt.
\newblock Practical and optimal {{LSH}} for angular distance.
\newblock In {\em Proceedings of the 29th {{International Conference}} on
  {{Neural Information Processing Systems}} - {{Volume}} 1}, volume~1 of {\em
  {{NIPS}}'15}, pages 1225--1233, Cambridge, MA, USA, December 2015. MIT Press.

\bibitem[AKK{\etalchar{+}}20]{agarwal2020leverage}
Naman Agarwal, Sham Kakade, Rahul Kidambi, Yin-Tat Lee, Praneeth Netrapalli,
  and Aaron Sidford.
\newblock Leverage score sampling for faster accelerated regression and erm.
\newblock In {\em Algorithmic Learning Theory}, pages 22--47. PMLR, 2020.

\bibitem[AL86]{axelsson1986rate}
Owe Axelsson and Gunhild Lindskog.
\newblock On the rate of convergence of the preconditioned conjugate gradient
  method.
\newblock {\em Numerische Mathematik}, 48:499--523, 1986.

\bibitem[BC13]{BurgisserCucker_Condition}
P.~B\"urgisser and F.~Cucker.
\newblock {\em Condition: The Geometry of Numerical Algorithms}.
\newblock Berlin, Heidelberg: Springer Nature, 2013.

\bibitem[Beb00]{bebendorf2000approximation}
Mario Bebendorf.
\newblock Approximation of boundary element matrices.
\newblock {\em Numer. Math.}, 86(4):565--589, 2000.

\bibitem[BEMOC21]{etayo-et-al}
Carlos Beltr{\'a}n, Uju{\'e} Etayo, Jordi Marzo, and Joaquim Ortega-Cerd{\`a}.
\newblock A sequence of polynomials with optimal condition number.
\newblock {\em J. Am. Math. Soc.}, 34(1):219--244, 2021.

\bibitem[BER04]{BER04}
Christopher Beattie, Mark Embree, and John Rossi.
\newblock Convergence of restarted {Krylov} subspaces to invariant subspaces.
\newblock {\em SIAM J.~Matrix Anal.\ Appl.}, 25:1074--1109, 2004.

\bibitem[BG65]{businger1965linear}
Peter Businger and Gene~H. Golub.
\newblock Linear least squares solutions by householder transformations.
\newblock {\em Numer. Math.}, 7(3):269–276, June 1965.

\bibitem[BGVKS23]{Banks-et-al-23}
Jess Banks, Jorge Garza-Vargas, Archit Kulkarni, and Nikhil Srivastava.
\newblock Pseudospectral shattering, the sign function, and diagonalization in
  nearly matrix multiplication time.
\newblock {\em Found. Comput. Math.}, 23(6):1959--2047, 2023.

\bibitem[BGVS22]{banks2022global}
Jess Banks, Jorge Garza-Vargas, and Nikhil Srivastava.
\newblock Global convergence of hessenberg shifted qr ii: Numerical stability.
\newblock {\em arXiv preprint arXiv:2205.06810}, 2022.

\bibitem[BHM00]{Briggs2000}
William~L. Briggs, Van~Emden Henson, and Steve~F. McCormick.
\newblock {\em A Multigrid Tutorial}.
\newblock Society for Industrial and Applied Mathematics, Philadelphia, PA, 2nd
  edition, 2000.

\bibitem[BLV18]{bai2018robust}
Zhaojun Bai, Ding Lu, and Bart Vandereycken.
\newblock Robust rayleigh quotient minimization and nonlinear eigenvalue
  problems.
\newblock {\em SIAM J.~Sci.\ Comput.}, 40(5):A3495--A3522, 2018.

\bibitem[BSS10]{BSS10}
M.~Bellalij, Y.~Saad, and H.~Sadok.
\newblock Futher analysis of the {Arnoldi} process for eigenvalue problems.
\newblock {\em SIAM J.~Numer.\ Anal.}, 48:393--407, 2010.

\bibitem[Car09]{Car09b}
Russell Carden.
\newblock A simple algorithm for the inverse field of values problem.
\newblock {\em Inverse Problems}, 25:115019 (9pp), 2009.

\bibitem[CC14]{chen2014stable}
Jie Chen and Edmond Chow.
\newblock A stable scaling of {Newton}-{Schulz} for improving the sign function
  computation of a {Hermitian} matrix.
\newblock {\em Preprint]. ANL/MCS-P5059-0114}, 2014.

\bibitem[CD13]{chiu2013sublinear}
Jiawei Chiu and Laurent Demanet.
\newblock Sublinear randomized algorithms for skeleton decompositions.
\newblock {\em SIAM J.~Matrix Anal.\ Appl.}, 34(3):1361--1383, 2013.

\bibitem[CDD25a]{chenakkod2025optimal}
Shabarish Chenakkod, Micha{\l} Derezi{\'n}ski, and Xiaoyu Dong.
\newblock Optimal oblivious subspace embeddings with near-optimal sparsity.
\newblock In {\em Proceedings of the 52nd International Colloquium on Automata,
  Languages, and Programming}, 2025.

\bibitem[CDD25b]{chenakkod2026optimal}
Shabarish Chenakkod, Micha{\l} Derezi{\'n}ski, and Xiaoyu Dong.
\newblock Optimal subspace embeddings: Resolving nelson-nguyen conjecture up to
  sub-polylogarithmic factors.
\newblock {\em arXiv preprint arXiv:2508.14234}, 2025.

\bibitem[CDDR24]{chenakkod2024optimal}
Shabarish Chenakkod, Micha{\l} Derezi{\'n}ski, Xiaoyu Dong, and Mark Rudelson.
\newblock Optimal embedding dimension for sparse subspace embeddings.
\newblock In {\em Proceedings of the 56th Annual ACM Symposium on Theory of
  Computing}, pages 1106--1117, 2024.

\bibitem[CE12]{CE12}
Russell~L. Carden and Mark Embree.
\newblock {Ritz} value localization for non-{Hermitian} matrices.
\newblock {\em SIAM J.~Matrix Anal.\ Appl.}, 33:1320--1338, 2012.

\bibitem[CEMT25]{CEMT25a}
Chris Cama{\~n}o, Ethan~N. Epperly, Raphael~A. Meyer, and Joel~A. Tropp.
\newblock Faster {{Linear Algebra Algorithms}} with {{Structured Random
  Matrices}}, August 2025.

\bibitem[CFG{\etalchar{+}}16]{CFG+16}
Krzysztof Choromanski, Francois Fagan, Cedric {Gouy-Pailler}, Anne Morvan,
  Tamas Sarlos, and Jamal Atif.
\newblock {{TripleSpin}} - a generic compact paradigm for fast machine learning
  computations, June 2016.

\bibitem[CH13]{CH13}
Russell Carden and Derek Hansen.
\newblock Ritz values of normal matrices and {Ceva}'s theorem.
\newblock {\em Linear Algebra Appl.}, 438:4114--4129, 2013.

\bibitem[Che24]{chen_24}
Tyler Chen.
\newblock The {L}anczos algorithm for matrix functions: a handbook for
  scientists, 2024.

\bibitem[CJHS25]{TN-sign-problem}
Jielun Chen, Jiaqing Jiang, Dominik Hangleiter, and Norbert Schuch.
\newblock Sign problem in tensor-network contraction.
\newblock {\em PRX Quantum}, 6:010312, Jan 2025.

\bibitem[CKM20]{cortinovis2020maximum}
Alice Cortinovis, Daniel Kressner, and Stefano Massei.
\newblock On maximum volume submatrices and cross approximation for symmetric
  semidefinite and diagonally dominant matrices.
\newblock {\em Linear Algebra Appl.}, 593:251--268, 2020.

\bibitem[CKS25]{ceruti2025low}
Gianluca Ceruti, Daniel Kressner, and Dominik Sulz.
\newblock Low-rank tree tensor network operators for long-range pairwise
  interactions.
\newblock {\em SIAM Journal on Scientific Computing}, 47(4):A2248--A2271, 2025.

\bibitem[CLS24]{CarLieStr24}
Erin Carson, J\"org Liesen, and Zden{\v e}k Strako{\v s}.
\newblock Towards understanding {CG} and {GMRES} through examples.
\newblock {\em Linear Algebra Appl.}, 692:241--291, 2024.

\bibitem[CNR{\etalchar{+}}25]{CNR+25}
Tyler Chen, Pradeep Niroula, Archan Ray, Pragna Subrahmanya, Marco Pistoia, and
  Niraj Kumar.
\newblock {{GPU-Parallelizable Randomized Sketch-and-Precondition}} for
  {{Linear Regression}} using {{Sparse Sign Sketches}}, June 2025.

\bibitem[Coh16]{cohen2016nearly}
Michael~B Cohen.
\newblock Nearly tight oblivious subspace embeddings by trace inequalities.
\newblock In {\em Proc. of the 27th annual ACM-SIAM Symposium on Discrete
  Algorithms}, pages 278--287. SIAM, 2016.

\bibitem[CT24]{chen_trogdon_24}
Tyler Chen and Thomas Trogdon.
\newblock Stability of the lanczos algorithm on matrices with regular spectral
  distributions.
\newblock {\em Linear Algebra and its Applications}, 682:191–237, February
  2024.

\bibitem[CY25]{cortinovis2025sublinear}
Alice Cortinovis and Lexing Ying.
\newblock A sublinear-time randomized algorithm for column and row subset
  selection based on strong rank-revealing {QR} factorizations.
\newblock {\em SIAM J.~Matrix Anal.\ Appl.}, 46(1):22--44, 2025.

\bibitem[Dav08]{davies2008approximate}
E.~Brian Davies.
\newblock Approximate diagonalization.
\newblock {\em SIAM J.~Matrix Anal.\ Appl.}, 29(4):1051--1064, 2008.

\bibitem[DDS24a]{Demmel_deflating_subspace}
James Demmel, Ioana Dumitriu, and Ryan Schneider.
\newblock Fast and inverse-free algorithms for deflating subspaces.
\newblock arXiv:2310.00193, 2024.

\bibitem[DDS24b]{Demmel2024}
James Demmel, Ioana Dumitriu, and Ryan Schneider.
\newblock Generalized pseudospectral shattering and inverse-free matrix pencil
  diagonalization.
\newblock {\em Foundations of Computational Mathematics}, 2024.

\bibitem[DG16]{drmac2016new}
Zlatko Drmac and Serkan Gugercin.
\newblock A new selection operator for the discrete empirical interpolation
  method---improved a priori error bound and extensions.
\newblock {\em SIAM J.~Sci.\ Comput.}, 38(2):A631--A648, 2016.

\bibitem[DGK98]{DGK98}
V.~Druskin, A.~Greenbaum, and L.~Knizhnerman.
\newblock Using nonorthogonal {L}anczos vectors in the computation of matrix
  functions.
\newblock {\em SIAM J. Sci. Comput.}, 19:38--54, 1998.

\bibitem[DGTY25]{damle2025estimating}
Anil Damle, Silke Glas, Alex Townsend, and Annan Yu.
\newblock Estimating a matrix's singular values with interpolative
  decompositions.
\newblock {\em Linear Algebra and its Applications}, 2025.

\bibitem[DH11]{DavisHu2011}
Timothy~A. Davis and Yifan Hu.
\newblock The university of florida sparse matrix collection.
\newblock {\em ACM Trans.\ Math.\ Software}, 38(1), December 2011.

\bibitem[Dhi97]{Dhillon:CSD-97-971}
Inderjit~Singh Dhillon.
\newblock {\em A New $O(n^2)$ Algorithm for the Symmetric Tridiagonal
  Eigenvalue/Eigenvector Problem}.
\newblock PhD thesis, EECS Department, University of California, Berkeley, Oct
  1997.

\bibitem[DK91]{druskin_knizhnerman_91}
V.~L. Druskin and L.~A. Knizhnerman.
\newblock Error bounds in the simple lanczos procedure for computing functions
  of symmetric matrices and eigenvalues.
\newblock {\em Comput. Math. Math. Phys.}, 31(7):20–30, 7 1991.

\bibitem[DLDM21]{derezinski2021sparse}
Michal Derezinski, Zhenyu Liao, Edgar Dobriban, and Michael Mahoney.
\newblock Sparse sketches with small inversion bias.
\newblock In {\em Conference on Learning Theory}, pages 1467--1510. PMLR, 2021.

\bibitem[DLNR25]{derezinski2025fine}
Michal Derezi{\'n}ski, Daniel LeJeune, Deanna Needell, and Elizaveta Rebrova.
\newblock Fine-grained analysis and faster algorithms for iteratively solving
  linear systems.
\newblock {\em Journal of Machine Learning Research}, 26(144):1--49, 2025.

\bibitem[DM23]{DM23}
Yijun Dong and Per-Gunnar Martinsson.
\newblock Simpler is better: A comparative study of randomized pivoting
  algorithms for {{CUR}} and interpolative decompositions.
\newblock {\em Advances in Computational Mathematics}, 49(4):66, August 2023.

\bibitem[DMMS11]{DMMS11}
Petros Drineas, Michael~W. Mahoney, S.~Muthukrishnan, and Tam{\'a}s Sarl{\'o}s.
\newblock Faster least squares approximation.
\newblock {\em Numerische Mathematik}, 117(2):219--249, February 2011.

\bibitem[DMSW05]{de2005near}
Stefano De~Marchi, Robert Schaback, and Holger Wendland.
\newblock Near-optimal data-independent point locations for radial basis
  function interpolation.
\newblock {\em Advances in Computational Mathematics}, 23(3):317--330, 2005.

\bibitem[DMY25]{derezinski2024faster}
Micha{\l} Derezi{\'n}ski, Christopher Musco, and Jiaming Yang.
\newblock Faster linear systems and matrix norm approximation via multi-level
  sketched preconditioning.
\newblock {\em ACM-SIAM Symposium on Discrete Algorithms (SODA)}, 2025.

\bibitem[DP03]{Dhillon_orthogonal}
Inderjit~S. Dhillon and Beresford~N. Parlett.
\newblock Orthogonal eigenvectors and relative gaps.
\newblock {\em SIAM J.~Matrix Anal.\ Appl.}, 25(3):858--899, 2003.

\bibitem[DP04]{DHILLON20041}
Inderjit~S. Dhillon and Beresford~N. Parlett.
\newblock Multiple representations to compute orthogonal eigenvectors of
  symmetric tridiagonal matrices.
\newblock {\em Linear Algebra and its Applications}, 387:1--28, 2004.

\bibitem[DPV04]{Dhillon:CSD-04-1346}
Inderjit Dhillon, Beresford~N. Parlett, and Christof Voemel.
\newblock Lapack working note 162: The design and implementation of the mrrr
  algorithm.
\newblock Technical Report UCB/CSD-04-1346, 2004.

\bibitem[DR24]{derezinski2024sharp}
Micha\l{} Derezi\'{n}ski and Elizaveta Rebrova.
\newblock {Sharp Analysis of Sketch-and-Project Methods via a Connection to
  Randomized Singular Value Decomposition}.
\newblock {\em SIAM J.~Math.\ Data Sci.}, 6(1):127--153, 2024.

\bibitem[DS99]{DanielsScuseria1999}
A.~D. Daniels and G.~E. Scuseria.
\newblock A recursive polynomial expansion method for the density matrix.
\newblock {\em Journal of Chemical Physics}, 110(3):1321--1328, 1999.

\bibitem[DS25]{derezinski2025approaching}
Micha{\l} Derezi{\'n}ski and Aaron Sidford.
\newblock Approaching optimality for solving dense linear systems with low-rank
  structure.
\newblock {\em arXiv preprint arXiv:2507.11724}, 2025.

\bibitem[DTM12]{DM12}
Jurjen Duintjer~Tebbens and G\'erard Meurant.
\newblock Any {Ritz} value behavior is possible for {Arnoldi} and {GMRES}.
\newblock {\em SIAM J.~Matrix Anal.\ Appl.}, 33:958--978, 2012.

\bibitem[DY24]{derezinski2024solving}
Micha{\l} Derezi{\'n}ski and Jiaming Yang.
\newblock Solving dense linear systems faster than via preconditioning.
\newblock In {\em 56th Annual ACM Symposium on Theory of Computing}, 2024.

\bibitem[EG11]{ErnstGander2011}
Oliver~G. Ernst and Martin~J. Gander.
\newblock Why it is difficult to solve {Helmholtz} problems with classical
  iterative methods.
\newblock In {\em Numerical Analysis of Multiscale Problems}, volume~83 of {\em
  Lecture Notes in Computational Science and Engineering}, pages 325--363.
  Springer, 2011.

\bibitem[EGLS21]{ehrlacher2021}
Virginie Ehrlacher, Laura Grigori, Damiano Lombardi, and Hao Song.
\newblock Adaptive hierarchical subtensor partitioning for tensor compression.
\newblock {\em SIAM Journal on Scientific Computing}, 43(1):A139--A163, 2021.

\bibitem[Emb22]{Emb22}
Mark Embree.
\newblock How descriptive are {GMRES} convergence bounds?
\newblock Technical Report arXiv2209:01231, 2022.

\bibitem[EMN25]{EMN25}
Ethan~N. Epperly, Maike Meier, and Yuji Nakatsukasa.
\newblock Fast randomized least-squares solvers can be just as accurate and
  stable as classical direct solvers.
\newblock {\em Communications on Pure and Applied Mathematics, to appear},
  2025.

\bibitem[Epp25]{Epp25a}
Ethan~N. Epperly.
\newblock {\em Make the Most of What You Have: {{Resource-efficient}}
  Randomized Algorithms for Matrix Computations}.
\newblock PhD thesis, California Institute of Technology, 2025.

\bibitem[ESW06]{ElmanHowardSilvester2014}
Howard Elman, David Silvester, and Andy Wathen.
\newblock {\em Finite Elements and Fast Iterative Solvers: With Applications in
  Incompressible Fluid Dynamics}.
\newblock Oxford University Press, 2014-06.

\bibitem[EVO04]{Erlangga2004}
Yogi~A. Erlangga, Cornelis Vuik, and Cornelis~W. Oosterlee.
\newblock On a class of preconditioners for solving the helmholtz equation.
\newblock {\em Applied Numerical Mathematics}, 50(3-4):409--425, 2004.

\bibitem[EY07]{EremenkoYuditskii2007}
A.~Eremenko and P.~Yuditskii.
\newblock Uniform approximation of sgn x by polynomials and entire functions.
\newblock {\em Journal of Analysis and Mathematics}, 101:313--324, 2007.

\bibitem[EY09]{engquist2009fast}
Bj{\"o}rn Engquist and Lexing Ying.
\newblock A fast directional algorithm for high frequency acoustic scattering
  in two dimensions.
\newblock 2009.

\bibitem[FG25a]{fahrbach2025tight}
Matthew Fahrbach and Mehrdad Ghadiri.
\newblock A tight lower bound for the approximation guarantee of higher-order
  singular value decomposition.
\newblock {\em arXiv preprint arXiv:2508.06693}, 2025.

\bibitem[FG25b]{FG25}
Israa Fakih and Laura Grigori.
\newblock Efficient {{QR-based Column Subset Selection}} through {{Randomized
  Sparse Embeddings}}, September 2025.

\bibitem[FGNW14]{Grigori2014}
Riadh Fezzani, Laura Grigori, Frédéric Nataf, and Ke~Wang.
\newblock Block filtering decomposition.
\newblock {\em Numerical Linear Algebra with Applications}, 21(6):703--721,
  2014.

\bibitem[FLa]{Fornace2025-approximation}
Mark Fornace and Michael Lindsey.
\newblock An approximation theory for {Markov} chain compression.

\bibitem[FLb]{Fornace2024-column-01}
Mark Fornace and Michael Lindsey.
\newblock Column and row subset selection using nuclear scores: algorithms and
  theory for nyst{r}öm approximation, {CUR} decomposition, and graph laplacian
  reduction.

\bibitem[FLT23]{FabLieTic23}
Vance Faber, J{\"o}rg Liesen, and Petr Tich{\'y}.
\newblock On the {F}orsythe conjecture.
\newblock {\em BIT}, 63:49, 2023.

\bibitem[For68]{For68}
George~E. Forsythe.
\newblock On the asymptotic directions of the {$s$}-dimensional optimum
  gradient method.
\newblock {\em Numer. Math.}, 11:57--76, 1968.

\bibitem[FY02]{Falgout2002}
Robert~D. Falgout and Ulrike~Meier Yang.
\newblock hypre: A library of high performance preconditioners.
\newblock In {\em Proceedings of the International Conference on Computational
  Science-Part III}, ICCS '02, page 632–641, Berlin, Heidelberg, 2002.
  Springer-Verlag.

\bibitem[GE96]{gu1996efficient}
Ming Gu and Stanley~C Eisenstat.
\newblock Efficient algorithms for computing a strong rank-revealing qr
  factorization.
\newblock {\em SIAM J.~Sci.\ Comput.}, 17(4):848--869, 1996.

\bibitem[GHRS18]{GowerHaRiSt2018}
Robert~M. Gower, Filip Hanzely, Peter Richt\'{a}rik, and Sebastian~U. Stich.
\newblock {Accelerated Stochastic Matrix Inversion: General Theory and Speeding
  up BFGS Rules for Faster Second-Order Optimization}.
\newblock In {\em Advances in Neural Information Processing Systems}, 2018.

\bibitem[GKS76]{golub1976rank}
G~Golub, V~Klema, and GW~Stewart.
\newblock Rank degeneracy and least squares problems.
\newblock In {\em Rep. STAN-CS-76-559}. Dept. Comput. Sci., Stanford
  University, 1976.

\bibitem[GL03]{GROER200345}
Benedikt Großer and Bruno Lang.
\newblock An $o(n^2)$ algorithm for the bidiagonal svd.
\newblock {\em Linear Algebra and its Applications}, 358(1):45--70, 2003.

\bibitem[GMNS19]{GerMadNieStr19}
Tom\'a\v~s{} Gergelits, Kent-Andr\'e{} Mardal, Bj{\o}rn~Fredrik Nielsen, and
  Zden{\v e}k Strako{\v s}.
\newblock Laplacian preconditioning of elliptic {PDE}s: localization of the
  eigenvalues of the discretized operator.
\newblock {\em SIAM J. Numer. Anal.}, 57(3):1369--1394, 2019.

\bibitem[GNS20]{GerNieStr20}
Tom\'a{\v s} Gergelits, Bj{\o}rn~Fredrik Nielsen, and Zden{\v e}k Strako{\v s}.
\newblock Generalized spectrum of second order differential operators.
\newblock {\em SIAM J. Numer. Anal.}, 58(4):2193--2211, 2020.

\bibitem[GNS22]{GerNieStr22}
Tom\'a\v~s{} Gergelits, Bj\o rn~Fredrik Nielsen, and Zden{\v e}k Strako{\v s}.
\newblock Numerical approximation of the spectrum of self-adjoint operators in
  operator preconditioning.
\newblock {\em Numer. Algorithms}, 91(1):301--325, 2022.

\bibitem[GR15]{GowerRichtarik2015}
Robert~M. Gower and Peter Richt\'{a}rik.
\newblock {Randomized Iterative Methods for Linear Systems}.
\newblock {\em SIAM J.~Matrix Anal.\ Appl.}, 36(4):1660--1690, 2015.

\bibitem[Gre89]{greenbaum_89}
Anne Greenbaum.
\newblock Behavior of slightly perturbed {L}anczos and conjugate-gradient
  recurrences.
\newblock {\em Linear Algebra and its Applications}, 113:7 -- 63, 1989.

\bibitem[Gre97]{Greenbaum97}
A.~Greenbaum.
\newblock Estimating the attainable accuracy of recursively computed residual
  methods.
\newblock {\em SIAM J. Matrix Anal. Appl.}, 18:535--551, 1997.

\bibitem[GS]{Guruswami2012-optimal}
Venkatesan Guruswami and Ali~Kemal Sinop.
\newblock {\em Optimal Column-Based Low-Rank Matrix Reconstruction}, pages
  1207--1214.
\newblock Proceedings. Society for Industrial and Applied Mathematics.

\bibitem[GS12]{GS12}
Thierry Gallay and Denis Serre.
\newblock Numerical measure of a complex matrix.
\newblock 65:287--336, 2012.

\bibitem[Hac85]{Hac85}
Wolfgang Hackbusch.
\newblock {\em Multigrid methods and applications}, volume~4 of {\em Springer
  Series in Computational Mathematics}.
\newblock Springer-Verlag, Berlin, 1985.

\bibitem[Hig94]{higham1994matrix}
Nicholas~J Higham.
\newblock The matrix sign decomposition and its relation to the polar
  decomposition.
\newblock {\em Linear Algebra and its Applications}, 212:3--20, 1994.

\bibitem[HJ25]{henning2025gross}
Patrick Henning and Elias Jarlebring.
\newblock The {Gross--Pitaevskii} equation and eigenvector nonlinearities:
  numerical methods and algorithms.
\newblock {\em SIAM Review}, 67(2):256--317, 2025.

\bibitem[HP92]{hong1992rank}
Yoo~Pyo Hong and C-T Pan.
\newblock Rank-revealing qr factorizations and the singular value
  decomposition.
\newblock {\em Math.\ Comp.}, 58(197):213--232, 1992.

\bibitem[HS52]{HestenesStiefel1952}
Magnus~R. Hestenes and Eduard Stiefel.
\newblock {Method of Conjugate Gradients for Solving Linear Systems}.
\newblock {\em Journal of Research of the National Bureau of Standards},
  49(6):409--436, 1952.

\bibitem[JCSH25]{TN-positive-bias}
Jiaqing Jiang, Jielun Chen, Norbert Schuch, and Dominik Hangleiter.
\newblock Positive bias makes tensor-network contraction tractable.
\newblock In {\em Proceedings of the 57th Annual ACM Symposium on Theory of
  Computing}, STOC '25, page 471–482, New York, NY, USA, 2025. Association
  for Computing Machinery.

\bibitem[JHNP13]{Jolivet2013}
Pierre Jolivet, Frédéric Hecht, Frédéric Nataf, and Christophe Prud'homme.
\newblock Scalable domain decomposition preconditioners for heterogeneous
  elliptic problems.
\newblock In {\em SC '13: Proceedings of the International Conference on High
  Performance Computing, Networking, Storage and Analysis}, pages 1--11, 2013.

\bibitem[JL25]{jarlebring2025polynomialsetassociatedfixed}
Elias Jarlebring and Gustaf Lorentzon.
\newblock The polynomial set associated with a fixed number of matrix-matrix
  multiplications, 2025.

\bibitem[JT25]{jeong2025convergence}
Sungwoo Jeong and Alex Townsend.
\newblock Convergence of pivoted {Cholesky} algorithm for lipschitz kernels.
\newblock {\em arXiv preprint arXiv:2509.13582}, 2025.

\bibitem[KB09]{kolda2009tensor}
Tamara~G Kolda and Brett~W Bader.
\newblock Tensor decompositions and applications.
\newblock {\em SIAM review}, 51(3):455--500, 2009.

\bibitem[KCP22]{dlr}
Jason Kaye, Kun Chen, and Olivier Parcollet.
\newblock Discrete lehmann representation of imaginary time green's functions.
\newblock {\em Phys. Rev. B}, 105:235115, Jun 2022.

\bibitem[KKR06]{kempe2006complexity}
Julia Kempe, Alexei Kitaev, and Oded Regev.
\newblock The complexity of the local {Hamiltonian} problem.
\newblock {\em SIAM J.~Comput.}, 35(5):1070--1097, 2006.

\bibitem[KN12]{KN12}
Daniel~M. Kane and Jelani Nelson.
\newblock Sparser {{Johnson-Lindenstrauss}} transforms.
\newblock In {\em Proceedings of the Twenty-Third Annual {{ACM-SIAM}} Symposium
  on {{Discrete}} Algorithms}, pages 1195--1206, 2012.

\bibitem[Kni96]{knizhnerman_96}
L.~A. Knizhnerman.
\newblock The simple {L}anczos procedure: Estimates of the error of the {G}auss
  quadrature formula and their applications.
\newblock {\em Comput. Math. Math. Phys.}, 36(11):1481–1492, 1 1996.

\bibitem[KS24]{kressner2024randomized}
Daniel Kressner and Nian Shao.
\newblock A randomized small-block {Lanczos} method for large-scale null space
  computations.
\newblock {\em arXiv preprint arXiv:2407.04634}, 2024.

\bibitem[KT24]{KT24}
Anastasia Kireeva and Joel~A. Tropp.
\newblock Randomized matrix computations: {{Themes}} and variations.
\newblock {\em 2023 CIME Summer School on Machine Learning , Cetraro, Italy,
  3-7 July 2023}, February 2024.

\bibitem[LL10]{LeventhalLewis2010}
Dennis Leventhal and Adrian~S. Lewis.
\newblock {Randomized Methods for Linear Constraints: Convergence Rates and
  Conditioning}.
\newblock {\em Mathematics of Operations Research}, 35(3):641--654, 2010.

\bibitem[LR25]{LokRebrova2025}
Jackie Lok and Elizaveta Rebrova.
\newblock {Subspace-constrained randomized coordinate descent for linear
  systems with good low-rank matrix approximations}.
\newblock 2025.

\bibitem[LS05]{LieStr05}
J.~Liesen and Z.~Strako{\v s}.
\newblock G{MRES} convergence analysis for a convection-diffusion model
  problem.
\newblock {\em SIAM J. Sci. Comput.}, 26(6):1989--2009, 2005.

\bibitem[LS13a]{lee2013efficient}
Yin~Tat Lee and Aaron Sidford.
\newblock Efficient accelerated coordinate descent methods and faster
  algorithms for solving linear systems.
\newblock In {\em 2013 ieee 54th annual symposium on foundations of computer
  science}, pages 147--156. IEEE, 2013.

\bibitem[LS13b]{LieStr13}
J\"org Liesen and Zden{\v e}k Strako{\v s}.
\newblock {\em Krylov Subspace Methods. Principles and Analysis}.
\newblock Numerical Mathematics and Scientific Computation. Oxford University
  Press, Oxford, 2013.

\bibitem[Mal92a]{Malyshev1}
A.~N. Malyshev.
\newblock Guaranteed accuracy in spectral problems of linear algebra. {I}.
\newblock {\em Siberian Advances in Mathematics}, 2(1):144--197, 1992.

\bibitem[Mal92b]{Malyshev2}
A.~N. Malyshev.
\newblock Guaranteed accuracy in spectral problems of linear algebra. {II}.
\newblock {\em Siberian Advances in Mathematics}, 2(2):153--204, 1992.

\bibitem[Mal04]{Mal04}
S.~M. Malamud.
\newblock Inverse spectral problem for normal matrices and the {Gauss--Lucas}
  theorem.
\newblock {\em Trans.\ Amer.\ Math.\ Soc.}, 357:4043--4064, 2004.

\bibitem[MBM{\etalchar{+}}23]{MBM+23}
Maksim Melnichenko, Oleg Balabanov, Riley Murray, James Demmel, Michael~W.
  Mahoney, and Piotr Luszczek.
\newblock {{CholeskyQR}} with {{Randomization}} and {{Pivoting}} for {{Tall
  Matrices}} ({{CQRRPT}}), November 2023.

\bibitem[McW60]{McWeeny1960}
R.~McWeeny.
\newblock Some recent advances in density matrix theory.
\newblock {\em Reviews of Modern Physics}, 32(4):335--369, 1960.

\bibitem[MDV20]{Marques_SVD}
Osni Marques, James Demmel, and Paulo~B. Vasconcelos.
\newblock Bidiagonal svd computation via an associated tridiagonal
  eigenproblem.
\newblock {\em ACM Trans. Math. Softw.}, 46(2), May 2020.

\bibitem[MG03]{miranian2003strong}
L~Miranian and Ming Gu.
\newblock Strong rank revealing {LU} factorizations.
\newblock {\em Linear algebra and its applications}, 367:1--16, 2003.

\bibitem[Min96]{minami1996local}
Nariyuki Minami.
\newblock Local fluctuation of the spectrum of a multidimensional anderson
  tight binding model.
\newblock {\em Communications in mathematical physics}, 177(3):709--725, 1996.

\bibitem[MMA18]{MMS2018}
Ca. Musco, Ch. Musco, and Sidford A.
\newblock Stability of the {L}anczos method for matrix function approximation,
  2018.

\bibitem[MMM24]{meyer2024unreasonable}
Raphael Meyer, Cameron Musco, and Christopher Musco.
\newblock On the unreasonable effectiveness of single vector krylov methods for
  low-rank approximation.
\newblock In {\em Proceedings of the 2024 Annual ACM-SIAM Symposium on Discrete
  Algorithms (SODA)}, pages 811--845. SIAM, 2024.

\bibitem[MNS{\etalchar{+}}18]{musco2018spectrum}
Cameron Musco, Praneeth Netrapalli, Aaron Sidford, Shashanka Ubaru, and David~P
  Woodruff.
\newblock Spectrum approximation beyond fast matrix multiplication: Algorithms
  and hardness.
\newblock In {\em 9th Innovations in Theoretical Computer Science Conference
  (ITCS 2018)}. Schloss-Dagstuhl-Leibniz Zentrum f{\"u}r Informatik, 2018.

\bibitem[MRK22]{masseirobolkressner2022}
Stefano Massei, Leonardo Robol, and Daniel Kressner.
\newblock Hierarchical adaptive low-rank format with applications to
  discretized partial differential equations.
\newblock {\em Numerical Linear Algebra with Applications}, 29(6):e2448, 2022.

\bibitem[Nes12]{Nesterov2012}
Yu~Nesterov.
\newblock {Efficiency of Coordinate Descent Methods on Huge-Scale Optimization
  Problems}.
\newblock {\em SIAM J.~Optim.}, 22(2):341--362, 2012.

\bibitem[Nie22]{nie2022matrix}
Zipei Nie.
\newblock Matrix anti-concentration inequalities with applications.
\newblock In {\em Proceedings of the 54th Annual ACM SIGACT Symposium on Theory
  of Computing}, pages 568--581, 2022.

\bibitem[Nik02]{Niklasson2002}
A.~M.~N. Niklasson.
\newblock Expansion algorithm for the density matrix.
\newblock {\em Physical Review B}, 66:155115, 2002.

\bibitem[NN13]{NN13a}
Jelani Nelson and Huy~L. Nguy{\^e}n.
\newblock {{OSNAP}}: {{Faster}} numerical linear algebra algorithms via sparser
  subspace embeddings.
\newblock In {\em 2013 Ieee 54th Annual Symposium on Foundations of Computer
  Science}, pages 117--126. IEEE, 2013.

\bibitem[NN14]{NN14}
Jelani Nelson and Huy~L. Nguyen.
\newblock Lower {{Bounds}} for {{Oblivious Subspace Embeddings}}.
\newblock In David Hutchison, Takeo Kanade, Josef Kittler, Jon~M. Kleinberg,
  Alfred Kobsa, Friedemann Mattern, John~C. Mitchell, Moni Naor, Oscar
  Nierstrasz, C.~Pandu~Rangan, Bernhard Steffen, Demetri Terzopoulos, Doug
  Tygar, Gerhard Weikum, Javier Esparza, Pierre Fraigniaud, Thore Husfeldt, and
  Elias Koutsoupias, editors, {\em Automata, {{Languages}}, and
  {{Programming}}}, volume 8572, pages 883--894. Springer Berlin Heidelberg,
  Berlin, Heidelberg, 2014.

\bibitem[NS24]{NieStr24}
Bj{\o}rn~Fredrik Nielsen and Zden{\v e}k Strako{\v s}.
\newblock A simple formula for the generalized spectrum of second order
  self-adjoint differential operators.
\newblock {\em SIAM Rev.}, 66(1):125--146, 2024.
\newblock Revised reprint of ``Generalized spectrum of second order
  differential operators'' [4128499].

\bibitem[Ose11]{oseledets2011tensor}
Ivan~V Oseledets.
\newblock Tensor-train decomposition.
\newblock {\em SIAM J.~Sci.\ Comput.}, 33(5):2295--2317, 2011.

\bibitem[Osi25]{Osinsky2025close}
A.I. Osinsky.
\newblock Close to optimal column approximation using a single svd.
\newblock {\em Linear Algebra Appl.}, 725:359--377, 2025.

\bibitem[Pai76]{paige_76}
Christopher~Conway Paige.
\newblock {Error Analysis of the {L}anczos Algorithm for Tridiagonalizing a
  Symmetric Matrix}.
\newblock {\em IMA Journal of Applied Mathematics}, 18(3):341--349, 12 1976.

\bibitem[Pai80]{paige_80}
Christopher~Conway Paige.
\newblock Accuracy and effectiveness of the {L}anczos algorithm for the
  symmetric eigenproblem.
\newblock {\em Linear Algebra and its Applications}, 34:235 -- 258, 1980.

\bibitem[PHY20]{pang2020interpolative}
Qiyuan Pang, Kenneth~L. Ho, and Haizhao Yang.
\newblock Interpolative decomposition butterfly factorization.
\newblock {\em sisc}, 42(2):A1097--A1115, 2020.

\bibitem[PS73]{doi:10.1137/0202007}
Michael~S. Paterson and Larry~J. Stockmeyer.
\newblock On the number of nonscalar multiplications necessary to evaluate
  polynomials.
\newblock {\em SIAM J.~Comput.}, 2(1):60--66, 1973.

\bibitem[PS08]{PS08}
Beresford Parlett and Gilbert Strang.
\newblock Matrices with prescribed {Ritz} values.
\newblock {\em Linear Algebra Appl.}, 428:1725--1739, 2008.

\bibitem[PV05]{Parlett:CSD-04-1367}
Beresford~N. Parlett and Christof Voemel.
\newblock Lapack working note 163: How the mrrr algorithm can fail on tight
  eigenvalue clusters.
\newblock Technical Report UCB/CSD-04-1367, EECS Department, University of
  California, Berkeley, Mar 2005.

\bibitem[PV21]{peng2021solving}
Richard Peng and Santosh Vempala.
\newblock Solving sparse linear systems faster than matrix multiplication.
\newblock In {\em Proceedings of the 2021 ACM-SIAM Symposium on Discrete
  Algorithms (SODA)}, pages 504--521. SIAM, 2021.

\bibitem[RT94]{RedTre94}
Satish~C. Reddy and Lloyd~N. Trefethen.
\newblock Pseudospectra of the convection-diffusion operator.
\newblock {\em SIAM J. Appl. Math.}, 54(6):1634--1649, 1994.

\bibitem[Rub11]{Rubensson2011}
E.~H. Rubensson.
\newblock Nonmonotonic recursive polynomial expansions for linear-scaling
  density-matrix construction.
\newblock {\em Journal of Chemical Theory and Computation}, 7(5):1233--1236,
  2011.

\bibitem[Saa03]{Saad2003}
Yousef Saad.
\newblock {\em {Iterative Methods for Sparse Linear Systems}}.
\newblock Society for Industrial and Applied Mathematics, Philadelphia, 2
  edition, 2003.

\bibitem[Saa11]{Saa11}
Yousef Saad.
\newblock {\em Numerical Methods for Large Eigenvalue Problems}.
\newblock SIAM, Philadelphia, second edition, 2011.

\bibitem[Sar06]{Sar06}
Tam{\'a}s Sarl{\'o}s.
\newblock Improved {{Approximation Algorithms}} for {{Large Matrices}} via
  {{Random Projections}}.
\newblock In {\em 2006 47th {{Annual IEEE Symposium}} on {{Foundations}} of
  {{Computer Science}} ({{FOCS}}'06)}, pages 143--152, October 2006.

\bibitem[Sch24]{Schneider_Thesis}
Ryan Schneider.
\newblock {\em Pseudospectral Divide-and-Conquer for the Generalized Eigenvalue
  Problem}.
\newblock PhD thesis, UC San Diego, 2024.

\bibitem[SER16]{ifiss2016}
David Silvester, Howard Elman, and Alison Ramage.
\newblock {I}ncompressible {F}low and {I}terative {S}olver {S}oftware ({IFISS})
  version 3.5, September 2016.
\newblock {\tt http://www.manchester.ac.uk/ifiss/}.

\bibitem[SH18]{santin6convergence}
Gabriele Santin and Bernard Haasdonk.
\newblock Convergence rate of the data-independent {P}-greedy algorithm in
  kernel-based approximation.
\newblock {\em Dolomites Research Notes on Approximation}, 6:83--100, 2018.

\bibitem[Sha25a]{shah2025hermitian}
Rikhav Shah.
\newblock Hermitian diagonalization in linear precision.
\newblock In {\em Proceedings of the 2025 Annual ACM-SIAM Symposium on Discrete
  Algorithms (SODA)}, pages 5599--5615. SIAM, 2025.

\bibitem[Sha25b]{Sha25}
Rikhav Shah.
\newblock The pseudospectrum of random compressions of matrices.
\newblock {\em arXiv preprint arXiv:2501.01418}, 2025.

\bibitem[Sob24]{sobczyk2024deterministic}
Aleksandros Sobczyk.
\newblock Deterministic complexity analysis of hermitian eigenproblems.
\newblock {\em arXiv preprint arXiv:2410.21550}, 2024.

\bibitem[Sor92]{Sor92}
D.~C. Sorensen.
\newblock Implicit application of polynomial filters in a $k$-step {Arnoldi}
  method.
\newblock {\em SIAM J.~Matrix Anal.\ Appl.}, 13:357--385, 1992.

\bibitem[SSZ24]{sparsepseudospectralshattering}
Rikhav Shah, Nikhil Srivastava, and Edward Zeng.
\newblock {Sparse Pseudospectral Shattering}, 2024.

\bibitem[Ste02]{Ste02}
G.~W. Stewart.
\newblock Backward error bounds for approximate {Krylov} subspaces.
\newblock {\em Linear Algebra Appl.}, 340:81--86, 2002.

\bibitem[TOS01]{Trottenberg2001}
Ulrich Trottenberg, Cornelis~W. Oosterlee, and Anton Schüller.
\newblock {\em Multigrid}, volume~33 of {\em Texts in Applied Mathematics}.
\newblock Academic Press, San Diego, 2001.
\newblock With contributions by A. Brandt, P. Oswald and K. Stüben.

\bibitem[Tow14]{townsend2014computing}
Alex Townsend.
\newblock {\em Computing with functions in two dimensions}.
\newblock PhD thesis, University of Oxford, 2014.

\bibitem[Tre13]{Trefethen2013}
Lloyd~N. Trefethen.
\newblock {\em Approximation Theory and Approximation Practice}.
\newblock SIAM, 1st edition, 2013.

\bibitem[Tro11]{Tro11}
Joel~A. Tropp.
\newblock Improved analysis of the subsampled randomized {{Hadamard}}
  transform.
\newblock {\em Advances in Adaptive Data Analysis}, 03(01n02):115--126, April
  2011.

\bibitem[Tro25]{Tro25}
Joel~A. Tropp.
\newblock Comparison theorems for the minimum eigenvalue of a random
  positive-semidefinite matrix, January 2025.

\bibitem[TT15]{townsend2015continuous}
Alex Townsend and Lloyd~N Trefethen.
\newblock Continuous analogues of matrix factorizations.
\newblock {\em Proc. R. Soc. A.}, 471(2173):20140585, 2015.

\bibitem[TYUC19]{TYUC19}
Joel~A. Tropp, Alp Yurtsever, Madeleine Udell, and Volkan Cevher.
\newblock Streaming {{Low-Rank Matrix Approximation}} with an {{Application}}
  to {{Scientific Simulation}}.
\newblock {\em SIAM Journal on Scientific Computing}, 41(4):A2430--A2463,
  January 2019.

\bibitem[Wes92]{Wesseling1992}
Pieter Wesseling.
\newblock {\em An Introduction to Multigrid Methods}.
\newblock Pure and Applied Mathematics. John Wiley \& Sons, Chichester; New
  York, 1992.

\bibitem[WL12a]{willems_etna}
Paul~R. Willems and Bruno Lang.
\newblock The {MR}$^3$-{GK} algorithm for the bidiagonal {SVD}.
\newblock {\em Electron. Trans. Numer. Anal.}, 39:1--21, 2012.

\bibitem[WL12b]{Willems_factorization}
Paul~R. Willems and Bruno Lang.
\newblock Twisted factorizations and qd-type transformations for the {MR}$^3$
  algorithm---new representations and analysis.
\newblock {\em SIAM J.~Matrix Anal.\ Appl.}, 33(2):523--553, 2012.

\bibitem[WL13]{Willems_framework}
Paul~R. Willems and Bruno Lang.
\newblock A framework for the {MR}$^3$ algorithm: Theory and implementation.
\newblock {\em SIAM J.~Sci.\ Comput.}, 35(2):A740--A766, 2013.

\bibitem[WLV06]{Willems_lapack}
Paul~R. Willems, Bruno Lang, and Christof V\"{o}mel.
\newblock Computing the bidiagonal svd using multiple relatively robust
  representations.
\newblock {\em SIAM J.~Matrix Anal.\ Appl.}, 28(4):907--926, 2006.

\bibitem[Woo14]{Woo14b}
David~P. Woodruff.
\newblock Sketching as a {{Tool}} for {{Numerical Linear Algebra}}.
\newblock {\em Foundations and Trends{\textregistered} in Theoretical Computer
  Science}, 10(1-2):1--157, 2014.

\bibitem[WT23]{wang2023operatorlearninghyperbolicpartial}
Christopher Wang and Alex Townsend.
\newblock Operator learning for hyperbolic partial differential equations,
  2023.

\bibitem[XZ17]{Xu2017}
Jinchao Xu and Ludmil Zikatanov.
\newblock Algebraic multigrid methods.
\newblock {\em Acta Numerica}, 26:591--721, 2017.

\bibitem[YGM09]{yang2009convergence}
Chao Yang, Weiguo Gao, and Juan~C Meza.
\newblock On the convergence of the self-consistent field iteration for a class
  of nonlinear eigenvalue problems.
\newblock {\em SIAM J.~Matrix Anal.\ Appl.}, 30(4):1773--1788, 2009.

\bibitem[Yse93]{Yserentant1993}
Harry Yserentant.
\newblock Old and new convergence proofs for multigrid methods.
\newblock {\em Acta Numerica}, pages 285--326, 1993.

\bibitem[ZB83]{ZhuBon83}
P.~F. Zhuk and L.~N. Bondarenko.
\newblock A conjecture of {G}. {E}. {F}orsythe.
\newblock {\em Mat. Sb. (N.S.)}, 121(163)(4):435--453, 1983.

\bibitem[ZM22]{ZouMag22}
Qinmeng Zou and Fr\'{e}d\'{e}ric Magoul\`es.
\newblock Delayed gradient methods for symmetric and positive definite linear
  systems.
\newblock {\em SIAM Review}, 64(3):517--553, 2022.

\end{thebibliography}

\end{document}